\documentstyle[epsf,amstex]{article}

\def\rA{Al} 
\def\rAra{Ar} 
\def\rB{B1} 
\def\rBc{B2} 
\def\rBMa{BM1} 
\def\rBMb{BM2} 
\def\rBMc{BM3} 
\def\rBMd{BM4} 
\def\rBMe{BM5} 
\def\rBMf{BM6} 
\def\rC{C} 
\def\rCJKS{CJKS} %
\def\rFRR{FRR} 
\def\rGoldA{G1} 
\def\rGoldB{G2} 
\def\rGPV{GPV} %
\def\rJ{J} 
\def\rJKS{JKS} 
\def\rKK{KK} 
\def\rKane{Kan} 
\def\rKauD{Kau1} %
\def\rKauE{Kau2} %
\def\rKauA{Kau3} %
\def\rKauB{Kau4} %
\def\rKauC{Kau5} %
\def\rKSW{KSW} 
\def\rMa{Ma} 
\def\rMUY{MUY} 
\def\rMor{Mo} %
\def\rNelson{N} %
\def\rSa{Sat} 
\def\rSaw{Saw} 
\def\rSWc{SW1} %
\def\rSWa{SW2} %
\def\rSWb{SW3} %
\def\rSk{Sk} 
\def\rTr{Tr} 
\def\rTu{Tu} 
\def\rVo{V} 
\def\rYa{Y} 

\def\fgCrossings{
\begin{figure}[htb]
\begin{center}
\mbox{\epsfxsize=2.0in \epsfbox{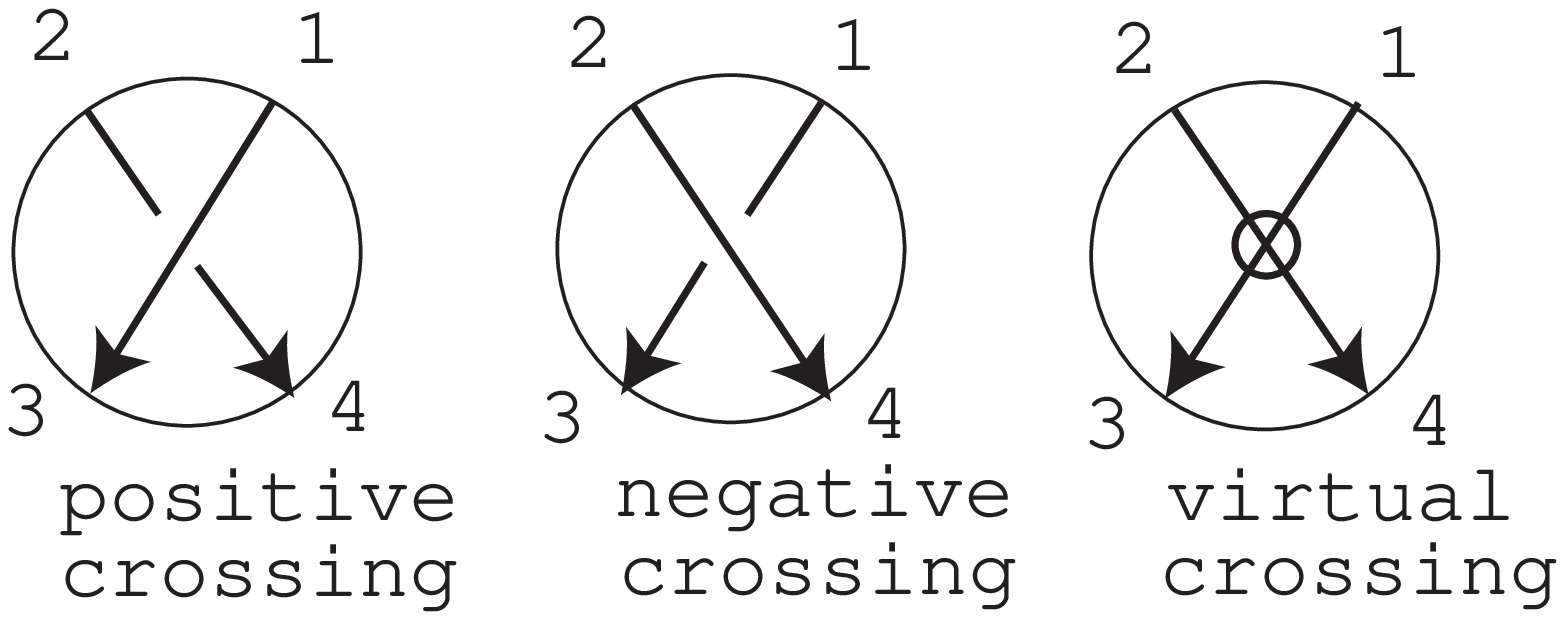}}
\end{center}
\caption{Crossings}
\label{fig:Crossings}\end{figure}}

\def\fgVBraidGene{
\begin{figure}[htb]
\begin{center}
\mbox{\epsfxsize=3.0in \epsfbox{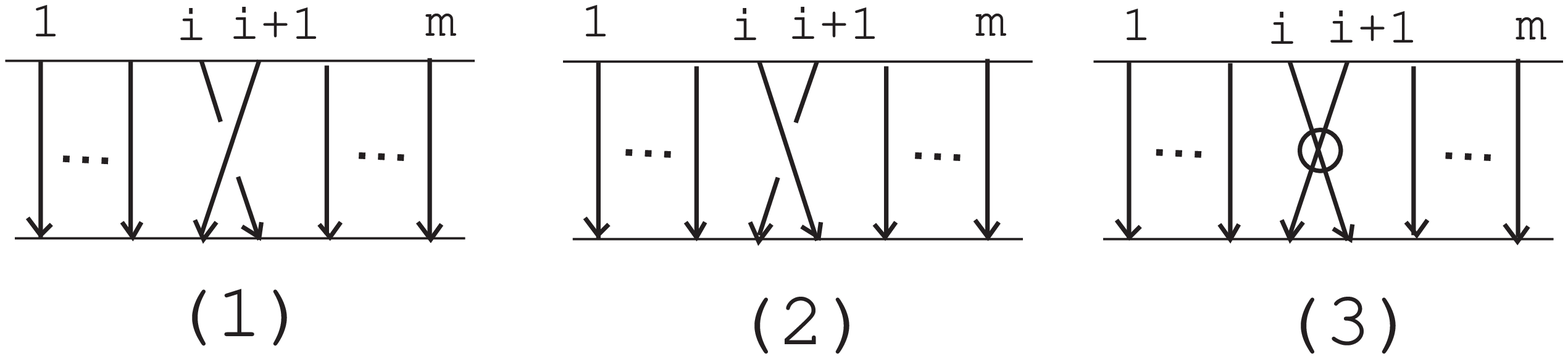}}
\end{center}
\caption{Standard generators}
\label{}\end{figure}}

\def\fgVRMove{
\begin{figure}[htb]
\begin{center}
\mbox{\epsfxsize=3.0in \epsfbox{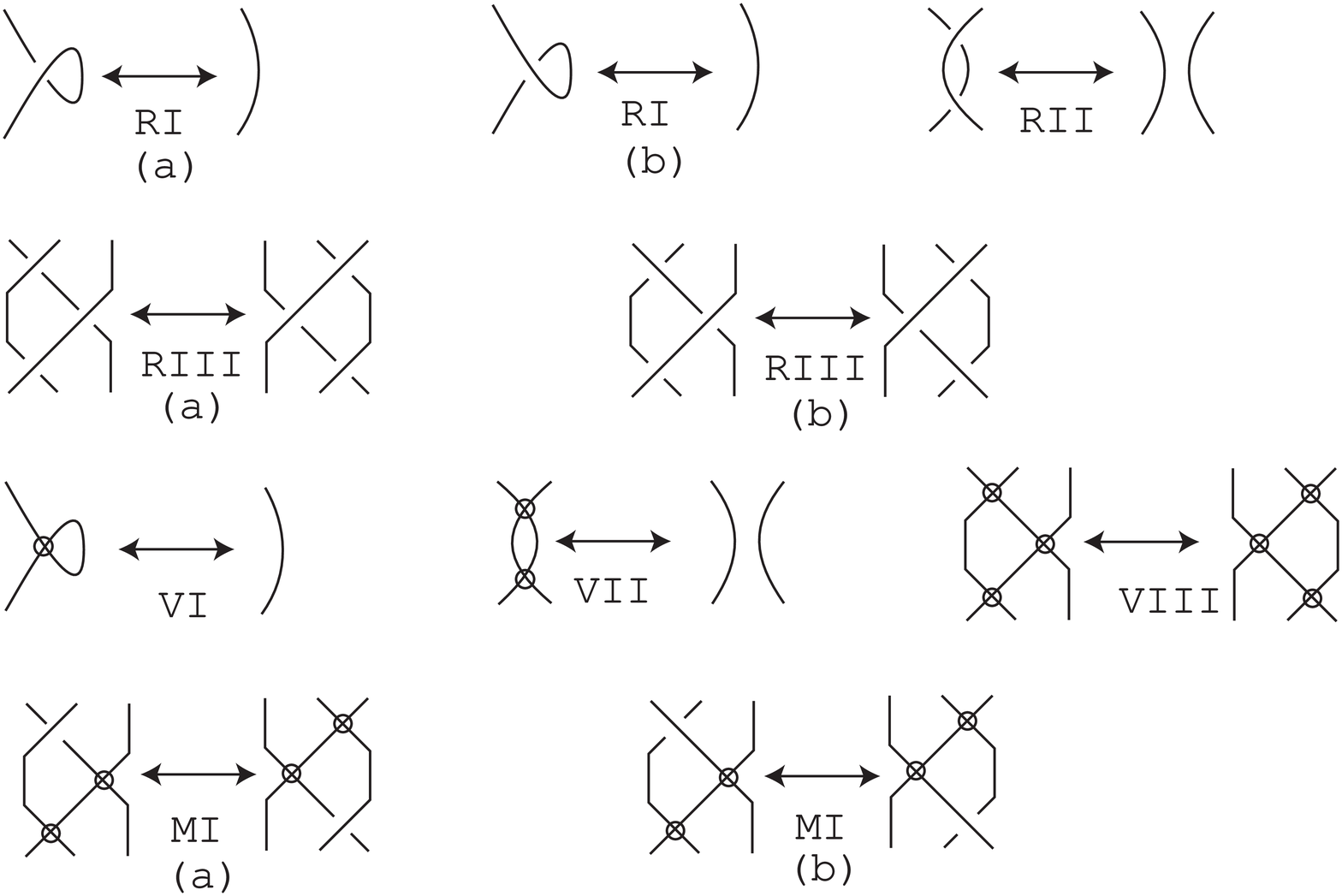}}
\end{center}
\caption{Virtual Reidemeister moves}
\label{}\end{figure}}

\def\fgVClosedBd{
\begin{figure}[htb]
\begin{center}
\mbox{\epsfxsize=2.5in \epsfbox{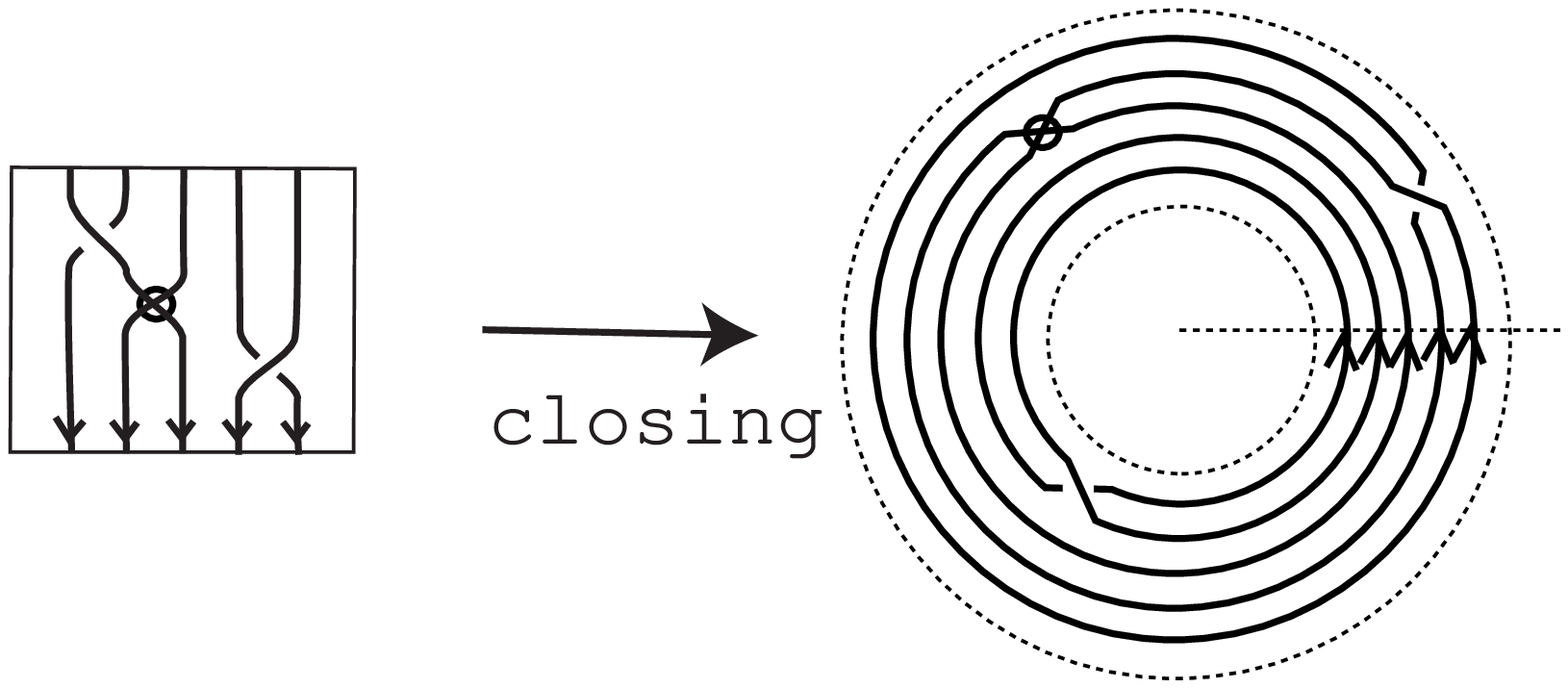}}
\end{center}
\caption{Closure}
\label{}\end{figure}}

\def\fgVStabliza{
\begin{figure}[htb]
\begin{center}
\mbox{\epsfxsize=3.0in \epsfbox{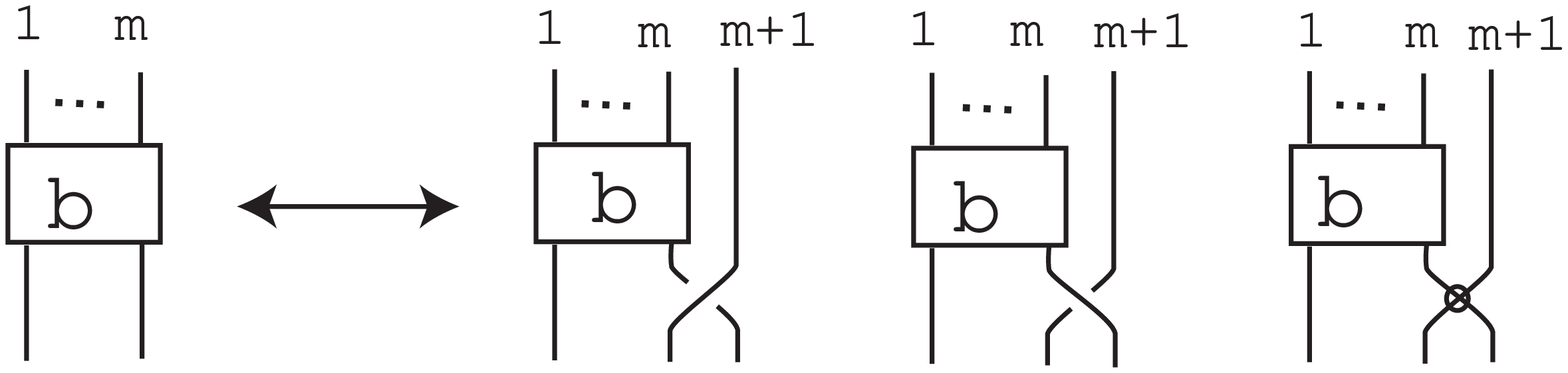}}
\end{center}
\caption{Right stabilizations}
\label{}\end{figure}}

\def\fgVExchangeM{
\begin{figure}[htb]
\begin{center}
\mbox{\epsfxsize=3.0in \epsfbox{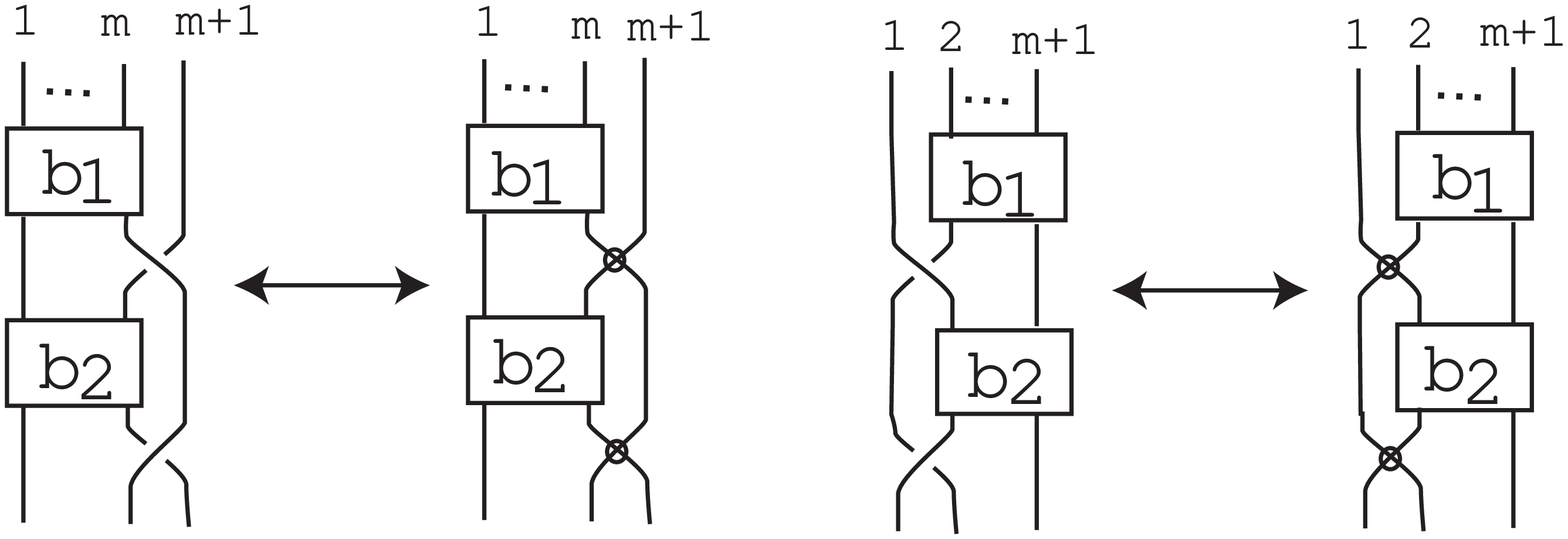}}
\end{center}
\caption{Right/left virtual exchange moves}
\label{}\end{figure}}

\def\fgVLinkD{
\begin{figure}[htb]
\begin{center}
\mbox{\epsfxsize=1.5in \epsfbox{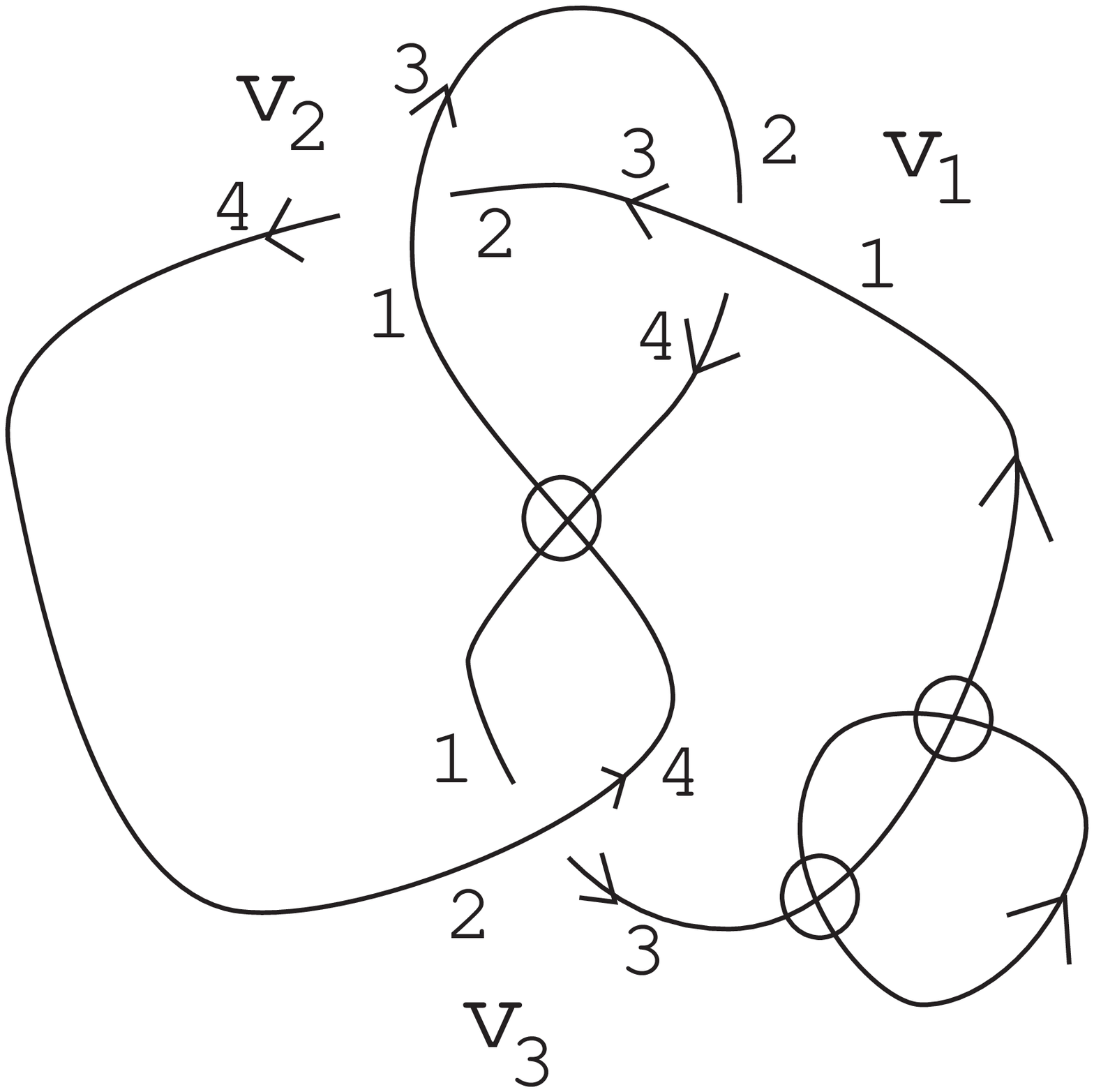}}
\end{center}
\caption{A virtual link diagram}
\label{}\end{figure}}

\def\fgarchomotopy{
\begin{figure}[htb]
\begin{center}
\mbox{\epsfxsize=3.5in \epsfbox{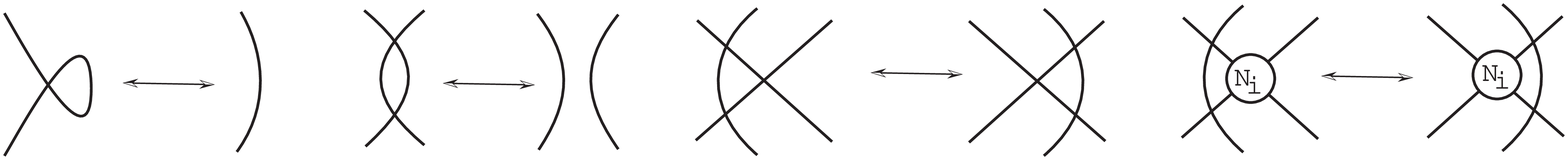}}
\end{center}
\caption{Moves on immersed curves}
\label{}\end{figure}}

\def\fgVStabB{
\begin{figure}[htb]
\begin{center}
\mbox{\epsfxsize=3.0in \epsfbox{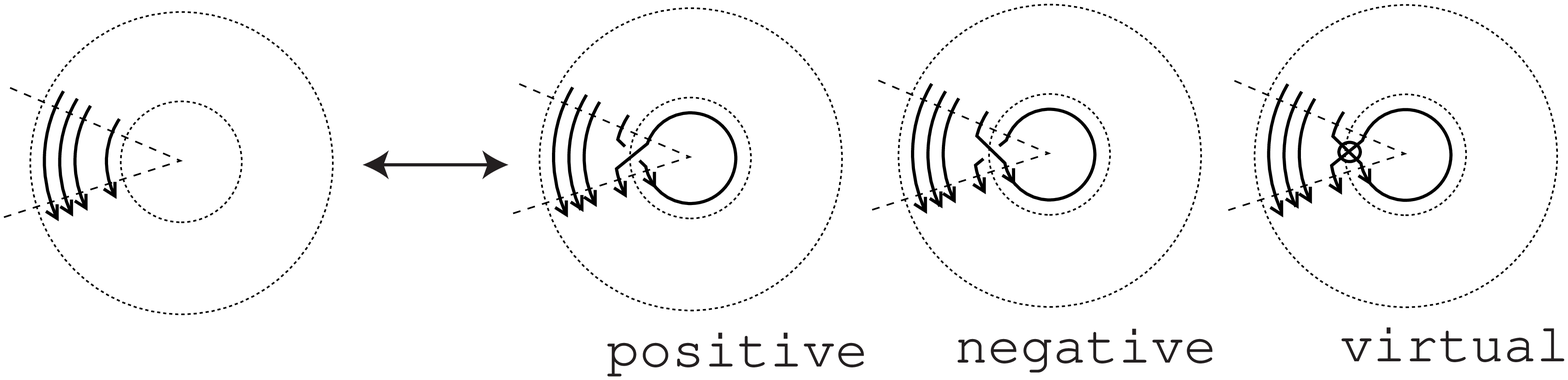}}
\end{center}
\caption{Right stabilizations (VM2-moves)}
\label{}\end{figure}}

\def\fgVBraidD{
\begin{figure}[htb]
\begin{center}
\mbox{\epsfxsize=4.0in \epsfbox{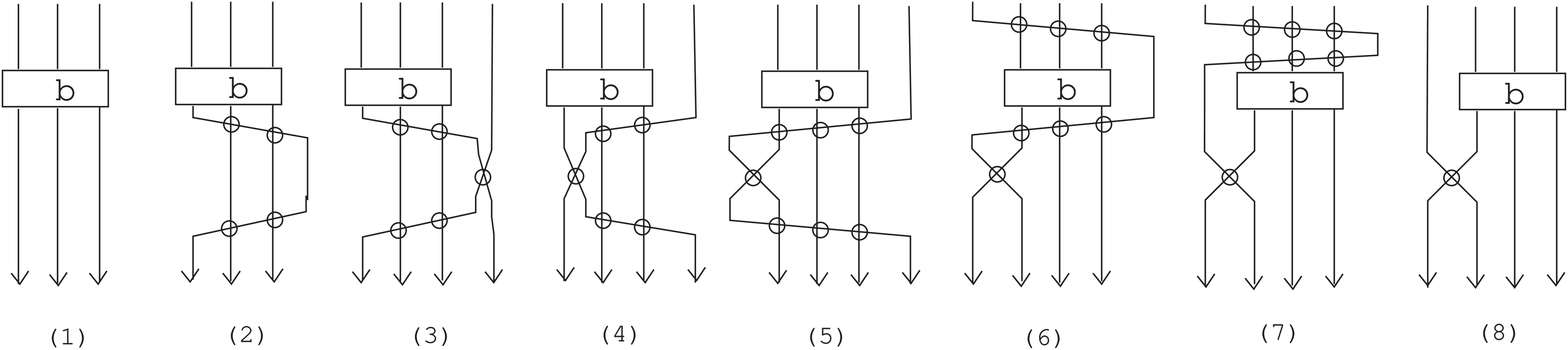}}
\end{center}
\caption{}
\label{fig:VBraidD}\end{figure}}

\def\fgVBraidA{
\begin{figure}[htb]
\begin{center}
\mbox{\epsfxsize=3.0in \epsfbox{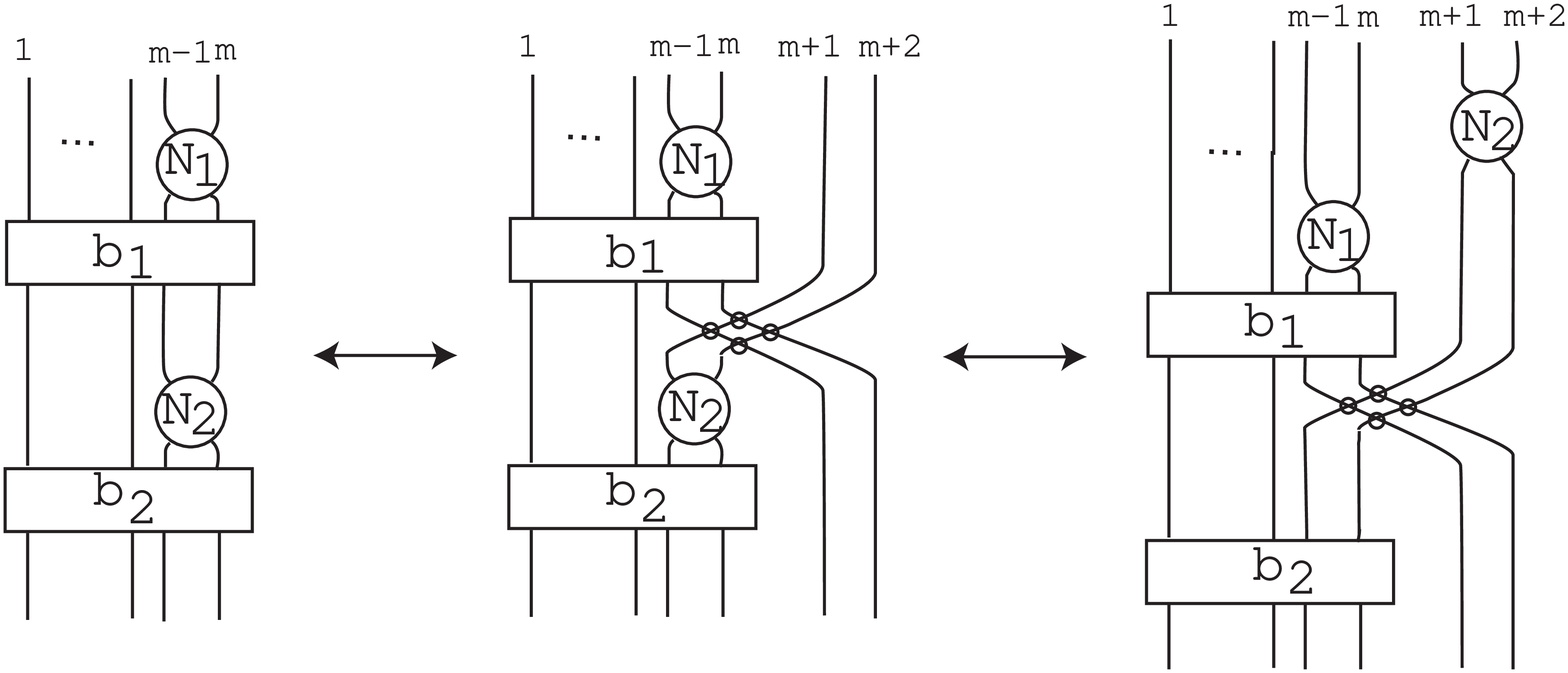}}
\end{center}
\caption{}
\label{fig:VBraidA}\end{figure}}

\def\fgVRMoveOriC{
\begin{figure}[htb]
\begin{center}
\mbox{\epsfxsize=3.0in \epsfbox{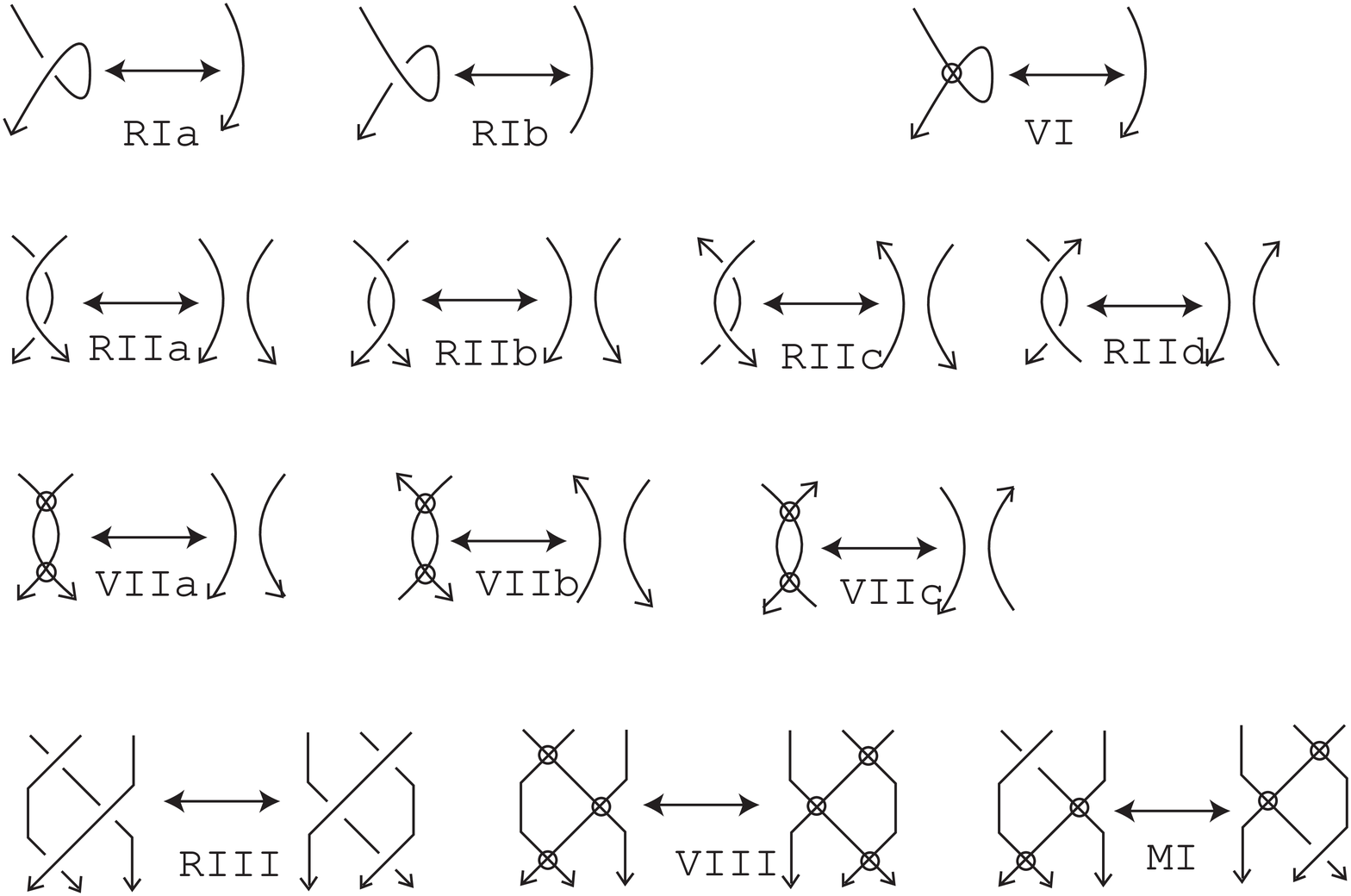}}
\end{center}
\caption{Oriented virtual Reidemeister moves}
\label{fig:VRMoveOriC}\end{figure}}

\def\fgVRMoveOriD{
\begin{figure}[htb]
\begin{center}
\mbox{\epsfxsize=2.0in \epsfbox{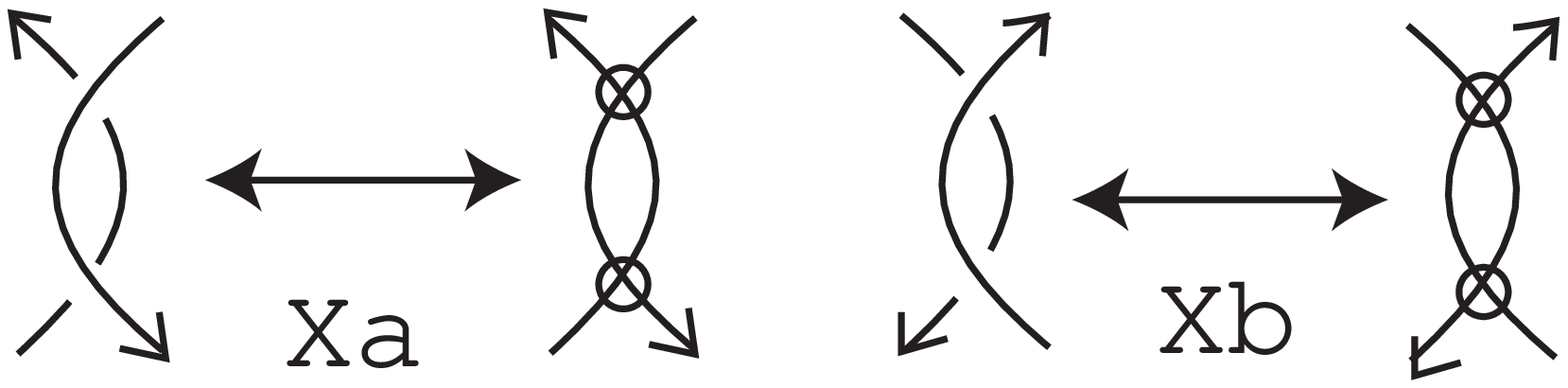}}
\end{center}
\caption{}
\label{fig:VRMoveOriD}\end{figure}}

\def\fgVBraidB{
\begin{figure}[htb]
\begin{center}
\mbox{\epsfxsize=3.0in \epsfbox{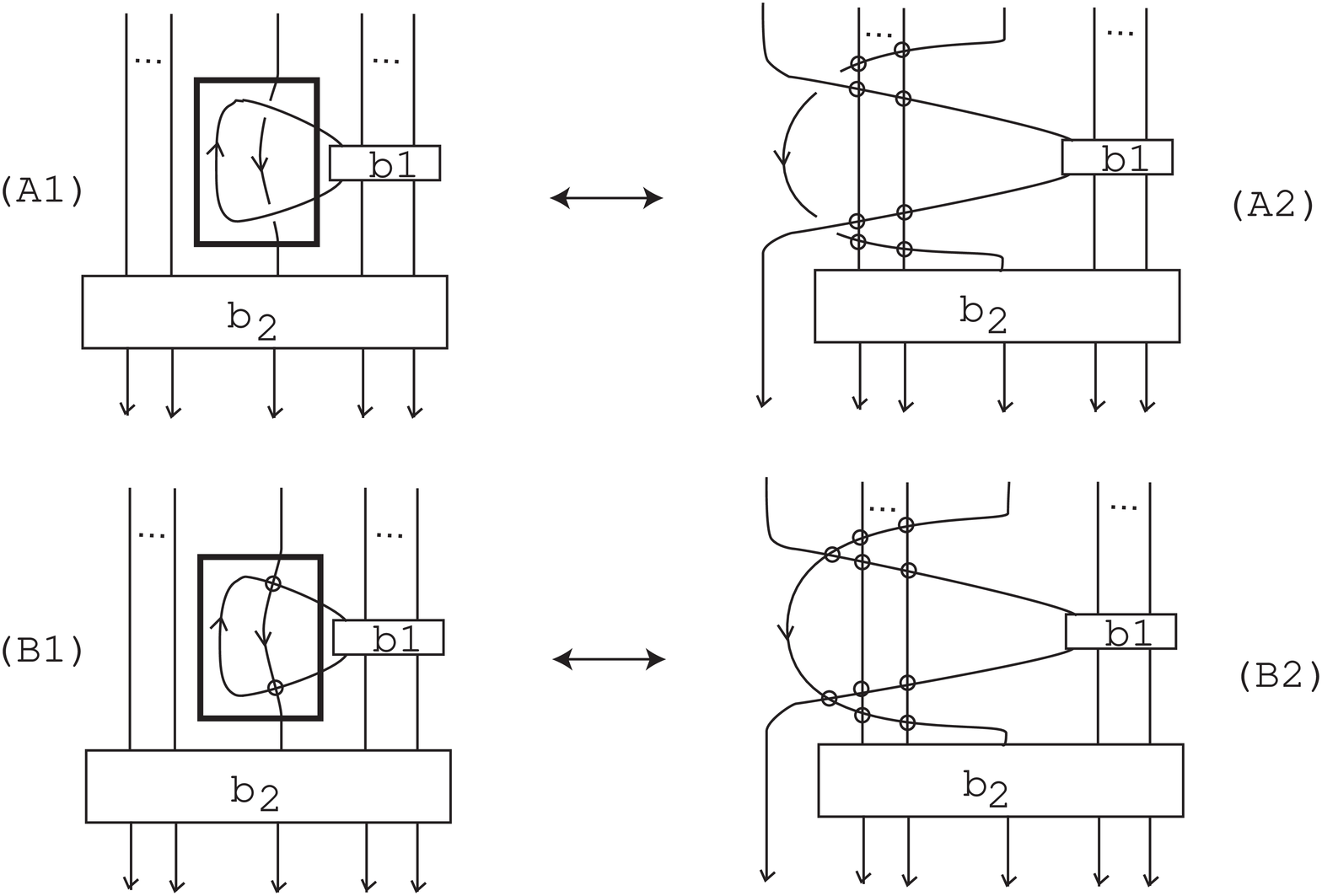}}
\end{center}
\caption{}
\label{fig:VBraidB}\end{figure}}

\def\fgVweights{
\begin{figure}[htb]
\begin{center}
\mbox{\epsfxsize=3.0in \epsfbox{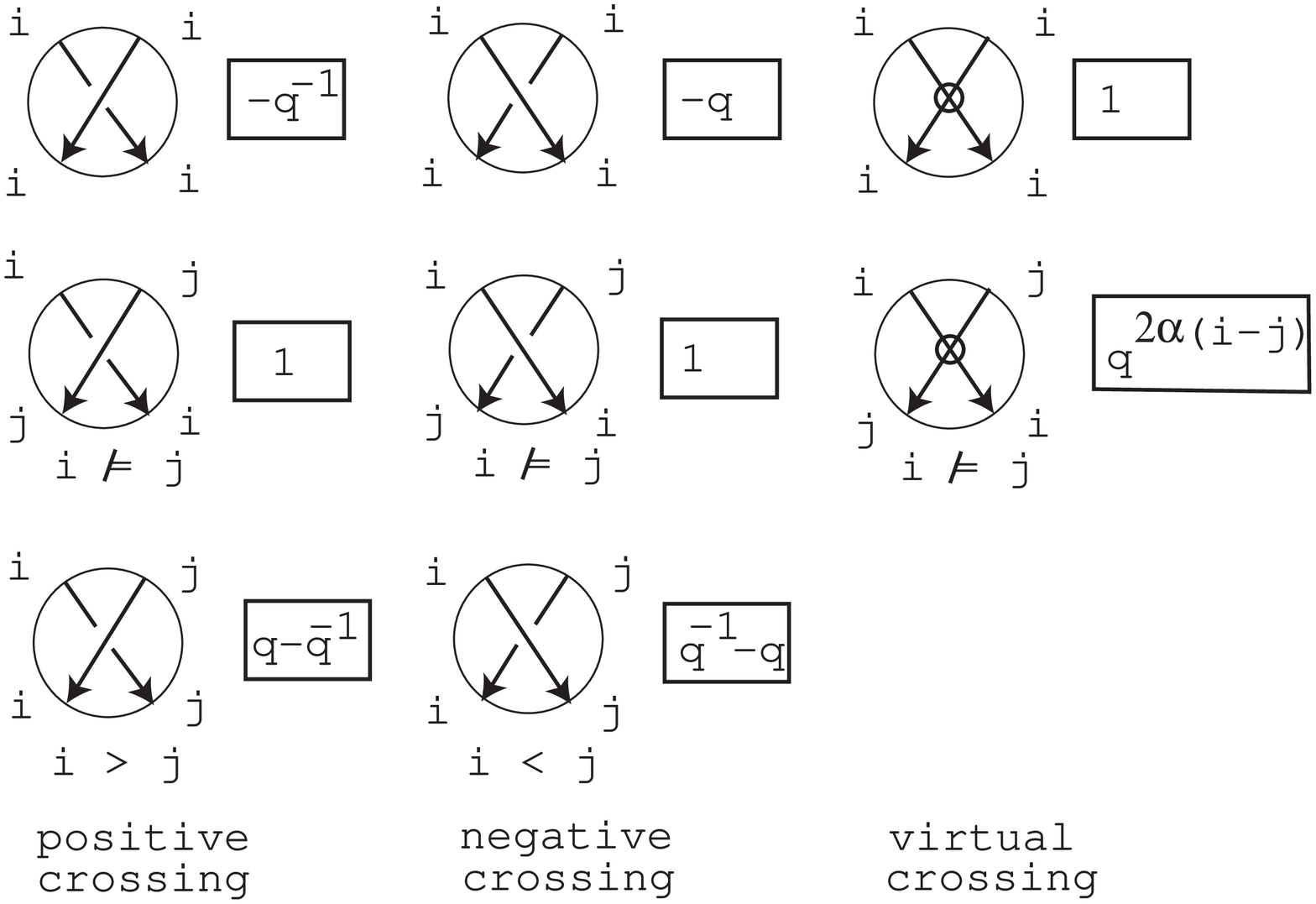}}
\end{center}
\caption{}
\label{fig:Vweights}\end{figure}}

\def\fgVsmooth{
\begin{figure}[htb]
\begin{center}
\mbox{\epsfxsize=2.5in \epsfbox{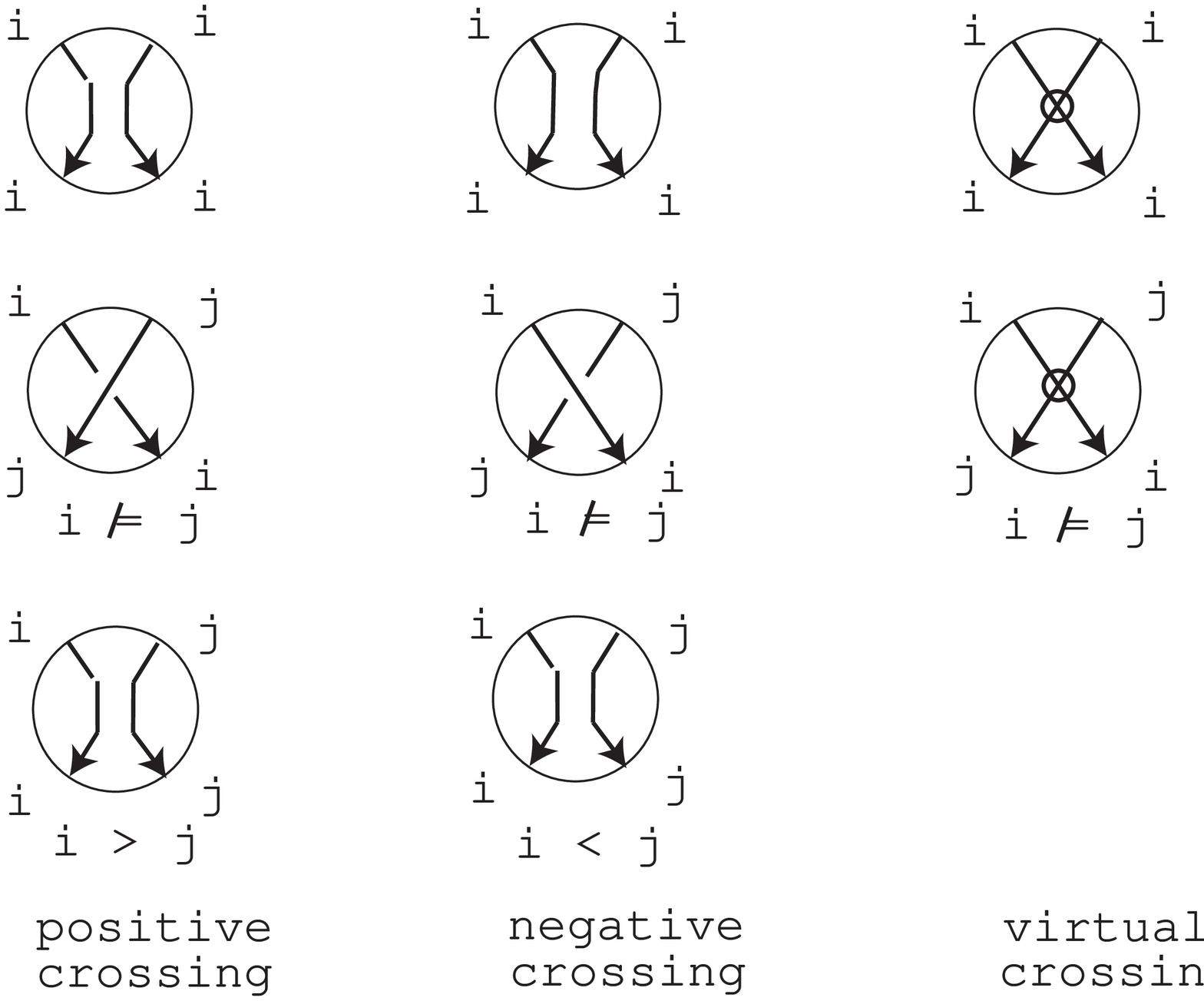}}
\end{center}
\caption{}
\label{fig:Vsmooth}\end{figure}}

\def\fgWRMoveA{
\begin{figure}[htb]
\begin{center}
\mbox{\epsfxsize=3.0in \epsfbox{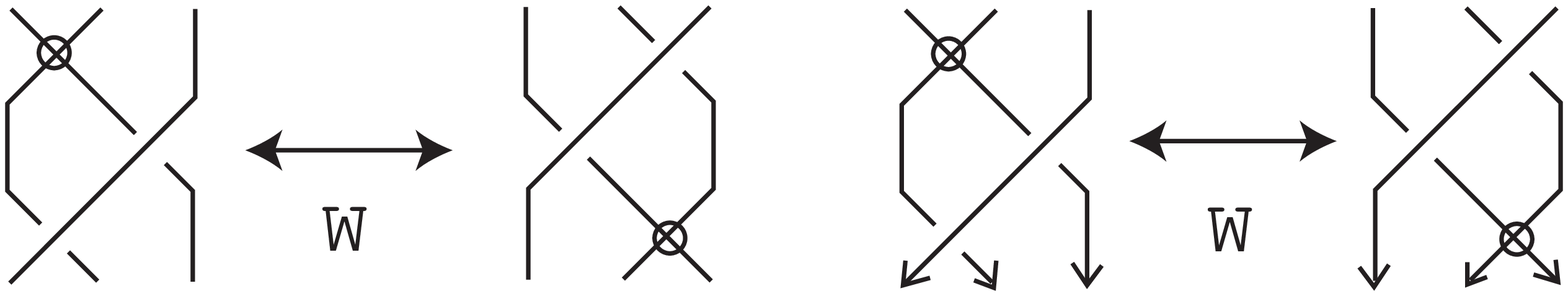}}
\end{center}
\caption{${\rm W}$-move}
\label{fig:WRMoveA}\end{figure}}

\def\fgWRMoveB{
\begin{figure}[htb]
\begin{center}
\mbox{\epsfxsize=3.0in \epsfbox{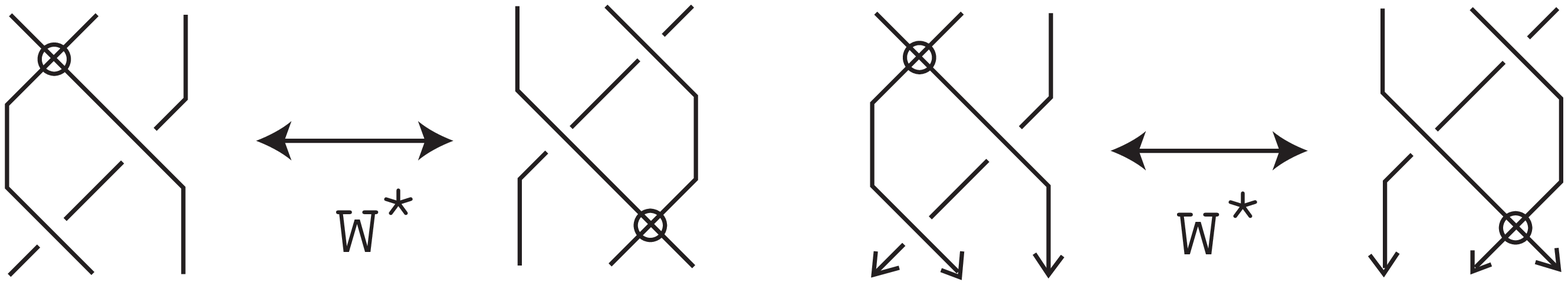}}
\end{center}
\caption{${\rm W}^*$-move}
\label{fig:WRMoveB}\end{figure}}

\def\qed{$\Box$}

\begin{document} 
\title{Braid presentation of virtual knots \\ and welded knots}
\author{Seiichi Kamada
\footnote{Supported by a Fellowship from 
the Japan Society for the Promotion of Science.} \\
Department of Mathematics, Osaka City University, \\
Sumiyoshi-ku, Osaka 558-8585, Japan}


\maketitle

\begin{abstract}
Virtual knots were first introduced by L. Kauffman, which are a  
generalization of classical knots and links. They lead us to 
the notion of virtual braids, which are closely related with welded braids 
of R. Fenn, R. Rim\'anyi and C. Rourke. 
It is proved that any virtual knot is uniquely described 
as the closure of a virtual braid modulo certain basic moves.  
This is analogous to the Alexander and Markov theorem 
for classical knots and braids.  
A similar result is proved for welded knots and braids
\end{abstract}

\noindent{\it Mathematics Subject Classification\/}:  
Primary 57M25. 

\noindent{\it Key words and phrases\/}: 
virtual knot, braid, welded braid, Alexander and Markov theorem, exchange move.

\section*{1. Introduction } 

In 1996, L.~Kauffman introduced 
the notion of a virtual knot, which is motivated by 
study of knots in a thickened surface and abstract Gauss codes, 
cf. \cite{\rKauA,\rKauB}. 
According to M.~Goussarov, M.~Polyak and O.~Viro \cite{\rGPV}, 
two classical knot diagrams represent the same knot type if and only if 
they represent the same virtual knot type. 
Thus, the notion of a virtual knot is a  
generalization of a classical knot in 3-space.  
Some properties and applications of virtual knots are found in 
\cite{\rCJKS, \rGPV, \rKK, \rKane, \rKauB, \rKSW, \rNelson, 
\rSa, \rSaw, \rSWc, \rSWa, \rSWb}, etc.  

Using the basic moves appearing in 
the definition of a virtual knot,  
we obtain the notion of a virtual braid 
(cf. \cite{\rKauB}).  
It is closely related with 
the welded braid group $WB_m$ and 
the braid-permutation group $BP_m$.  
R. Fenn, R. Rim\'anyi and C. Rourke \cite{\rFRR} defined   
the groups $WB_m$ and $BP_m$ and proved that they are isomorphic. 
There is a canonical epimorphism from the 
virtual braid group $VB_m$ to the welded braid group $WB_m$, and 
the group $VB_m$ contains  
the braid group $B_m$ and the symmetric group $S_m$ 
as subsets in a natural way. 

Braid theory plays an important role in classical knot theory.  
The two theories are related by 
Alexander's and Markov's theorems which state that every knot (or link) type is 
represented by the closure of a braid 
and such a braid presentation is unique up to 
conjugations and stabilizations (cf. \cite{\rA, \rB, \rBc, 
\rBMa, \rBMb, \rBMc, \rBMd, \rBMe, \rBMf, \rMa, \rSk, \rTr, \rVo, \rYa}).  
There is a one-to-one correspondence between 
links and braids modulo these operations.  
Analogously, virtual braid theory is expected to be so in virtual knot theory.  
It is quite easy to prove the Alexander theorem for virtual knots; that is, 
every virtual link type is represented by the closure of a 
virtual braid.  In fact, this is 
obvious from the relationship between virtual links and 
Gauss diagrams, or Gauss codes given in \cite{\rGPV,\rKauB}.  
In \cite{\rKauC}, Kauffman asked whether there is a generalization 
of the Markov theorem for virtual knots.  
Our main result is the following, which is an 
answer to his question and     
ensures a relationship between virtual braid theory and virtual knot theory.  

\noindent{\bf Theorem 3.2} {\it 
Two virtual braids (or virtual braid diagrams)  
have equivalent closures as virtual links  
if and only if they are related by a finite sequence of the 
following moves $({\rm VM1})$ -- $({\rm VM3})$ (or $({\rm VM0})$ -- $({\rm VM3})$): 
\begin{itemize}
\setlength{\itemsep}{-3pt}
\item[{\rm (VM0)}]  a braid move (which is a move  
corresponding to a defining relation of the virtual braid group), 
\item[{\rm (VM1)}] a conjugation (in the virtual braid group), 
\item[{\rm (VM2)}] a right stabilization of positive, negative or virtual type, 
and its inverse operation, 
\item[{\rm (VM3)}] a right/left virtual exchange move.   
\end{itemize} }

${\rm VM0}$-, ${\rm VM1}$- and ${\rm VM2}$-moves correspond to 
classical Markov moves.  
The last move (${\rm VM3}$-move) is an analogue 
of an exchange move (cf. \cite{\rBc, \rBMd}).  
In the category of classical braids,  an exchange move 
is a consequence of Markov moves. 
However its analogy does not hold 
in the category of virtual braids (Proposition~3.3).   
Thus ${\rm VM3}$-moves are essential.  
By Theorem~3.2, a left stabilization 
of any type (defined in \S~3) is realized by 
${\rm VM0}$-, ${\rm VM1}$-, ${\rm VM2}$- and ${\rm VM3}$-moves.  
If the left stabilization is of virtual type, 
then it can be realized  
without ${\rm VM3}$-moves (Proposition~3.4).    
This is analogous to a fact that 
a left stabilization of positive/negative 
type for classical braids is realized by Markov moves.  
It is rather surprising that a left stabilization of positive/negative 
type for virtual braids is not realized by ${\rm VM0}$-, ${\rm VM1}$- and ${\rm VM2}$-moves (Proposition~3.5).

Most of known virtual knot invariants, as groups, quandles, Alexander polynomials, 
and f-polynomials (Jones polynomials, normalized bracket polynomials), are 
considered and calculated easily via virtual braids 
(see \cite{\rKauB, \rSWa} for these invariants):   Group presentations and
quandle presentations are obtained from virtual braids  
by a method which is completely analogous to that in classical knot theory 
(although upper presentation 
and lower presentation yield different groups and quandles in general, 
\cite{\rGPV, \rKauB}).  
The Burau representation is easily defined for 
virtual braids and it brings Alexander polynomials.   
The notion of Temperley-Lieb algebra 
is generalized to virtual ones and its basis is much simpler 
than the classical one, \cite{\rKSW}.  
Using the virtual Temperley-Lieb algebra, 
one can obtain the f-polynomials via braids.  
Very recently, D. Silver and S. Williams \cite{\rSWb} 
found a new group invariant for 
a virtual link (with $\mu$ components) and 
a $(\mu +1)$-variable polynomial invariant derived from it.  
J.~Sawollek \cite{\rSaw} also found a similar invariant.  
Their invariants are so powerful as to distinguish the trivial knot 
and Kauffman's example \cite{\rKauB,\rKauC} of a virtual knot 
which cannot be distinguished 
by all of the above invariants. 
(J. S. Carter proved this fact independently 
by an argument in \cite{\rC}).  
D. Silver informed the author that their invariant in \cite{\rSWb} 
was motivated by the Burau representation and 
it is natural to consider via virtual braids.  
Related topics to these invariants from a point of view of virtual braids 
will be discussed elsewhere.  
                               
Welded braid theory due to Fenn, Rim\'anyi and Rourke \cite{\rFRR} 
yields the notion of 
welded knots and links, which are also known as 
virtual knots and links in the weak sense (\cite{\rKauC, \rSa}). 
We prove Alexander's and Markov's theorems for welded knots and links.   

\noindent{\bf Theorem 7.2} {\it 
Two welded braids (or welded braid diagrams)  
have equivalent closures as welded links  
if and only if they are related by a finite sequence of the 
following moves $({\rm WM1})$ -- $({\rm WM2})$ (or $({\rm WM0})$ -- $({\rm WM2})$): 
\begin{itemize} 
\setlength{\itemsep}{-3pt}
\item[{\rm (WM0)}] a welded braid move (which is a move 
corresponding to a defining relation of the welded braid group), 
\item[{\rm (WM1)}] a conjugation in the welded braid group, 
\item[{\rm (WM2)}] a right stabilization of positive, negative or virtual type, 
and its inverse operation. 
\end{itemize} }

The author wishes to thank J. S. Carter, N. Kamada, 
L. H. Kauffman, D. S. Silver, 
X.-S. Lin and O. Dasbach for many stimulating conversations.  
He also thanks to Department of Mathematics and Statistics, 
University of South Alabama for hospitality. 

\section*{2. Virtual braids }

Let $m$ be a positive integer and $Q_m$ a set of $m$ interior points 
of the interval $[0,1]$.  We denote by $E$ the 2-disk $[0,1]\times [0,1]$ 
and by $p_2:E \to [0,1]$ the second factor projection. 

\noindent{\bf Definition 2.1} 
A {\it virtual braid diagram of degree $m$\/} is an 
immersed 1-manifold $b = a_1 \cup \dots \cup a_m$ in $E$  
such that 
\begin{itemize} 
\setlength{\itemsep}{-3pt}
\item[{\rm (1)}] $\partial b= Q_m \times \{0,1\} \subset E$, 
\item[{\rm (2)}] for each $i \in \{1,\dots, m\}$, $p_2|_{a_i} : a_i \to [0,1]$ is a homeomorphism, 
\item[{\rm (3)}] the singularity (the multiple point set) $V(b)$ 
consists of transverse double points, 
\item[{\rm (4)}] $p_2|_{V(b)}: V(b) \to [0,1]$ is injective, 
\item[{\rm (5)}] each point of $V(b)$ is assigned information of 
{\it positive, negative \/} or {\it virtual crossing\/}  as in Figure~1 
(where the labels $1, \dots, 4$ are used later). 
\end{itemize}
The arcs $a_1, \dots, a_m$ are assumed to be oriented 
from the top ($[0,1]\times\{1\}$) to the bottom ($[0,1]\times\{0\}$) of $E$ and 
two virtual braid diagrams are identified if they are transformed  
into each other continuously keeping the above conditions.

\fgCrossings

The set of virtual braid diagrams of degree $m$ (with the concatenation product) 
forms a monoid which is 
generated by $\sigma_i, \sigma_i^{-1}, \tau_i$ ($i=1,\dots, m-1$) 
as in Figure~2.  The trivial element is $Q_m \times[0,1] \subset E$. 

\fgVBraidGene

\noindent{\bf Definition 2.2} 
The {\it virtual braid group $VB_m$ of degree $m$\/} is the group 
obtained from the monoid of virtual braid diagrams of degree $m$ 
by introducing the following relations: 
$$
\begin{array}{l}
{\rm (Trivial~relations)}  \quad
    \sigma_i \sigma_i^{-1} = \sigma_i^{-1} \sigma_i =1 \\
{\rm (Braid~relations)} \quad
  \cases
    \sigma_i \sigma_j = \sigma_j \sigma_i, \quad |i-j|>1 \\
    \sigma_i \sigma_{i+1} \sigma_i = \sigma_{i+1} \sigma_i \sigma_{i+1} 
  \endcases \\
{\rm (Permutation~group~relations)}  \quad
  \cases
    \tau_i^2 =1 \\
    \tau_i \tau_j = \tau_j \tau_i, \quad |i-j|>1 \\
    \tau_i \tau_{i+1} \tau_i = \tau_{i+1} \tau_i \tau_{i+1} 
  \endcases \\
{\rm (Mixed~relations)} \quad
  \cases
    \sigma_i \tau_j = \tau_j \sigma_i, \quad |i-j|>1 \\
    \sigma_i \tau_{i+1} \tau_i = \tau_{i+1} \tau_i \sigma_{i+1}.   
  \endcases
\end{array}
$$
A {\it  virtual braid of degree $m$\/} is an element of $VB_m$.  We 
denote it by the same symbol $b$ as its representative 
(a virtual braid diagram) $b$ unless it makes confusion.

The welded braid group $WB_m$ (defined in \cite{\rFRR}) is obtained 
from $VB_m$ by introducing additional relations  
$ \tau_i \sigma_{i+1} \sigma_i = \sigma_{i+1} \sigma_i \tau_{i+1} $ 
($i=1,\dots, m-2)$.  
There is a canonical epimorphism $VB_m \to WB_m$.  
In particular, we see that the subgroup of $VB_m$ generated by 
$\sigma_i$ ($i=1,\dots, m$) is isomorphic to the braid group $B_m$ 
and the subgroup generated by 
$\tau_i$ ($i=1,\dots, m$) is isomorphic to the symmetric group $S_m$, 
cf. \cite{\rFRR}.

\section*{3. Braid Presentation of Virtual Knots }

A {\it virtual link diagram\/} is a closed oriented 1-manifold 
$K$ immersed in $\text{\bf R}^2$ such that the singularity 
set $V(K)$ consists of transverse double points each of which 
is assigned information of positive, negative or virtual crossing 
as in Figure~1.  Positive and negative 
crossings are also called {\it real crossings\/}.  
Virtual link diagrams are considered up to isotopy of $\text{\bf R}^2$. 
{\it Virtual Reidemeister moves} are the 
local moves illustrated in Figure~3.  
Notice that the moves indicated by (b) are 
obtained from the moves indicated by (a) by use of   
RII-moves or VII-moves.   
Two virtual link diagrams are {\it equivalent\/} (as virtual links) if they are 
related by a finite sequence of virtual Reidemeister moves.   
A {\it virtual link\/} (or a {\it virtual link type\/}) is 
the equivalence class of a virtual link diagram, \cite{\rGPV, \rKauA, \rKauB}.   

\fgVRMove

For a virtual braid diagram $b$, the {\it closure\/} of $b$ is 
a virtual link diagram constructed as in Figure~4.    
If $b$ and $b'$ are equivalent as virtual braids, then 
their closures are equivalent as virtual links.  Thus 
in virtual knot theory the closure makes sense for a virtual braid.  

\fgVClosedBd

\noindent{\bf Proposition 3.1} {\it 
Every virtual link type is represented by 
the closure of a virtual braid diagram.  
}

For virtual braids $b_1, b_2 \in VB_m$,   
we say that $b_2^{-1} b_1 b_2$ is obtained from $b_1$ by a 
{\it conjugation\/}.  

For a virtual braid (diagram) $b$ of degree $m$, we denote by $\iota_s^t(b)$ 
the virtual braid (diagram) of degree $m+s+t$ obtained from $b$ by 
adding $s$ trivial arcs to the left of $b$ and $t$ 
trivial arcs to the right.  This defines a monomorphism 
$\iota_s^t: VB_m \to VB_{m+s+t}$.

A {\it right stabilization\/} 
of {\it positive}, {\it negative} or {\it virtual type\/} is  
a replacement of $b \in VB_m$ by
$\iota_0^1(b)\sigma_m$, 
$\iota_0^1(b)\sigma_m^{-1}$ or  
$\iota_0^1(b)\tau_m \in VB_{m+1}$, 
respectively.  
See Figure~5. 
Similarly, a {\it left stabilization\/} is 
a replacement of $b \in VB_m$ by 
$\iota_1^0(b)\sigma_1$,  
$\iota_1^0(b)\sigma_1^{-1}$ or 
$\iota_1^0(b)\tau_1$.  

\fgVStabliza

A {\it right virtual exchange move\/} is a replacement 
$$ \iota_0^1(b_1) \sigma_m^{-1} \iota_0^1(b_2) \sigma_m 
\quad \leftrightarrow \quad 
\iota_0^1(b_1) \tau_m \iota_0^1(b_2) \tau_m \quad \in VB_{m+1} $$
and  
a {\it left virtual exchange move\/} is a replacement 
$$ \iota_1^0(b_1) \sigma_1^{-1} \iota_1^0(b_2) \sigma_1
\quad \leftrightarrow \quad 
\iota_1^0(b_1) \tau_1 \iota_1^0(b_2) \tau_1 \quad \in VB_{m+1} $$
where $b, b' \in VB_m$, see Figure~6.

\fgVExchangeM

\noindent{\bf Theorem 3.2} {\it 
Two virtual braids (or virtual braid diagrams)  
have equivalent closures as virtual links  
if and only if they are related by a finite sequence of the 
following moves $({\rm VM1})$ -- $({\rm VM3})$ (or $({\rm VM0})$ -- $({\rm VM3})$): 
\begin{itemize} 
\setlength{\itemsep}{-3pt}
\item[{\rm (VM0)}] a braid move (which is a move 
corresponding to a defining relation of the virtual braid group), 
\item[{\rm (VM1)}] a conjugation (in the virtual braid group), 
\item[{\rm (VM2)}] a right stabilization of positive, negative or virtual type, 
and its inverse operation, 
\item[{\rm (VM3)}] a right/left virtual exchange move.   
\end{itemize} }

A replacement 
$$\iota_0^1(b_1) \sigma_m^{-1} \iota_0^1(b_2) \sigma_m 
\quad \leftrightarrow \quad
\iota_0^1(b_1) \sigma_m \iota_0^1(b_2) \sigma_m^{-1} \quad \in B_{m+1},$$ 
where $b_1, b_2 \in B_m$,   
is called an {\it exchange move}, cf. \cite{\rBc, \rBMd}. 
A ${\rm VM3}$-move is an analogue of this move.  
In classical braid theory, 
an exchange move is a consequence of braid moves, 
conjugations and right stabilizations.  
However, in virtual braid theory, 
a virtual exchange move is independent of these moves.  

\noindent{\bf Proposition 3.3} {\it 
A ${\rm VM3}$-move is not a consequence of 
${\rm VM0}$-,  ${\rm VM1}$- and ${\rm VM2}$-moves; 
namely, there is a pair of virtual braids 
which are related by a ${\rm VM3}$-move and 
never related by a sequence of ${\rm VM0}$-,  
${\rm VM1}$- and ${\rm VM2}$-moves.  
}

By Theorem~3.2, a left stabilization of any type for virtual braids is a 
consequence of ${\rm VM0}$-, ${\rm VM1}$-, ${\rm VM2}$- and ${\rm VM3}$-moves.  
If it is of virtual type, then we do not need ${\rm VM3}$-moves.  

\noindent{\bf Proposition 3.4} {\it 
A left stabilization of virtual type is a consequence of 
a ${\rm VM2}$-move and some ${\rm VM0}$- and ${\rm VM1}$-moves.  
}

This is analogous to a fact that 
a left stabilization (of positive/negative type) for classical braids is 
a consequence of a right stabilization and some braid moves and 
conjugations.  
If the left stabilization for virtual braids 
is of positive/negative type, then we need ${\rm VM3}$-moves in general.  

\noindent{\bf Proposition 3.5} {\it 
A left stabilization of positive/negative type 
for virtual braids 
is not a consequence of 
${\rm VM0}$-, ${\rm VM1}$- and ${\rm VM2}$-moves; 
namely, 
there is a pair of virtual braids 
which are related by a left stabilization 
of positive/negative type and 
never related by a sequence of ${\rm VM0}$-,  
${\rm VM1}$- and ${\rm VM2}$-moves.
}

\section*{4. Braiding Process}

For a virtual link diagram $K$, we denote by $V_R(K)$ the 
set of real crossings and by 
$S(K) : V_R(K) \to \{+1, -1\}$ the map assigning the real crossings 
their signs.  
For a real crossing 
$v \in V_R(K)$, let $N(v)$ be a regular neighborhood of $v$ 
as in Figure~1. 
We denote by $v^{(1)}, v^{(2)}, v^{(3)}, v^{(4)}$ 
the four points of $\partial N(v) \cap K$ ordered as in the figure.   
Put $W = W(K) = {\rm Cl}(\text{\bf R}^2 - \cup_{v \in V_R(K)} N(v))$ 
and $V_R^\partial (K) = \{ v^{(j)} \/ | \/ v \in V_R(K), j \in 
\{1,2,3,4\} \}$, where ${\rm Cl}$ means the closure. 
The restriction of $K$ to $W$, denoted by $K|_W$, is the union 
of some oriented arcs and loops 
immersed in $W$ such that the singularities are virtual crossings of $K$ and 
the boundaries of the arcs are the set $V_R^\partial (K)$.  

Define a subset $G(K) \subset    V_R^\partial (K) \times V_R^\partial (K)$ such that 
$(a,b) \in G(K)$ if and only if $K|_W$ has an arc starting 
from $a$ and terminating at $b$.  
We denote by $\mu(K)$ the number of components of $K$.
For example, for a virtual link diagram illustrated in Figure~7, 
$$
\begin{array}{cl}
V_R(K) &= \{ v_1, v_2, v_3 \},  \\
S(K) &: v_1 \mapsto +1, \quad v_2 \mapsto +1, \quad v_3 \mapsto -1, \\
G(K) &= \{ (v_3^{(3)}, v_1^{(1)}), (v_1^{(3)}, v_2^{(2)}), 
(v_2^{(4)}, v_3^{(2)}), (v_3^{(4)}, v_2^{(1)}), 
(v_2^{(3)}, v_1^{(2)}), (v_1^{(4)}, v_3^{(1)}) \},  \\
\mu(K) &= 2. 
\end{array}$$

\fgVLinkD

The {\it Gauss data\/} of $K$ is 
the quadruple $(V_R(K), S(K), G(K), \mu(K))$. 
We say that two virtual link diagrams $K$ and $K'$ have the 
{\it same Gauss data\/} if $\mu(K) = \mu(K')$ and if 
there is a bijection $g: V_R(K) \to V_R(K')$ such that 
$g$ preserves the signs of the crossing points and that 
$(a,b) \in G(K)$ implies $(g(a), g(b)) \in G(K')$, where 
$g: V_R^\partial (K) \to V_R^\partial (K')$ is the bijection 
induced from $g: V_R(K) \to V_R(K')$.  
This condition is equivalent to that 
$K$ and $K'$ have the same Gauss diagram in 
the sense of \cite{\rGPV} or the same Gauss code in the sense of 
\cite{\rKauB}. 

Let $K$ be a virtual link diagram and $W = W(K) = 
{\rm Cl}(\text{\bf R}^2 - \cup_{v \in V_R(K)} N(v))$ as before.  
Suppose that $K'$ is a 
virtual link diagram with the same Gauss data as $K$.  Then 
we can deform $K'$ by an isotopy of $\text{\bf R}^2$ 
such that 
\begin{itemize} 
\setlength{\itemsep}{-3pt}
\item[{\rm (1)}] $K$ and $K'$ are identical in $N(v)$ for every $v \in V_R(K)$, 
\item[{\rm (2)}] $K'$ has no real crossings in $W$, and 
\item[{\rm (3)}] there is a one-to-one correspondence between the arcs/loops 
of $K|_W$ and those of $K'|_W$ with respect to the 
end points of the arcs. 
\end{itemize}
In this situation, we say that $K'$ is obtained 
from $K$ by {\it replacing $K|_W$\/}.  

\noindent{\bf Lemma 4.1 {\rm (\cite{\rGPV, \rKauA, \rKauB})}} {\it 
If two virtual link diagrams $K$ and $K'$ 
have the same Gauss data, then $K$ is 
equivalent to $K'$.  
Moreover, such an equivalence can be realized by 
${\rm VI}$-, ${\rm VII}$-, ${\rm VIII}$- and ${\rm MI}$-moves. 
}

{\it Proof.} 
Without loss of generality we may assume that 
$K'$ is obtained from $K$ by replacing $K|_W$.  
Let $a_1, a_2, \dots, a_s$ be the arcs/loops of $K|_W$, 
and let $a'_1, a'_2, \dots, a'_s$ 
be corresponding ones for $K'|_W$.  We may assume that 
$a'_1$ intersects $a_2, \dots, a_s$ transversely. 
The arc or loop $a_1$ is homotopic to $a'_1$ in $\text{\bf R}^2$  
(relative to the boundary of $a_1$ if $a_1$ is an arc).  
Taking the homotopy generically 
with respect to the arcs/loops $a_2, \dots, a_s$ and the 2-disks 
$N_1, \dots, N_n$, we see that the arc/loop $a_1$ is 
transformed into 
$a'_1$ by a finite sequence of moves as in Figure~8  
up to isotopy of $\text{\bf R}^2$.   
Each move is a VI-, VII-, VIII-, or MI-move.  Inductively, 
every $a_i$ is transformed into $a'_i$ by such moves. 
\qed 

\fgarchomotopy

Let $O$ be the origin of $\text{\bf R}^2$.  
Identify $\text{\bf R}^2 -\{O\}$ with $\text{\bf R}_+ \times S^1$ 
by the polar coordinate and let $\pi: \text{\bf R}^2 -\{O\} = 
\text{\bf R}_+ \times S^1 \to S^1$ be the projection, where 
$\text{\bf R}_+$ is the half-line consisting of positive numbers.  
A {\it braided virtual link diagram} 
(of degree $m$) is a virtual link diagram $K$ such that 
\begin{itemize} 
\setlength{\itemsep}{-3pt}
\item[{\rm (1)}] it is contained in $\text{\bf R}^2 -\{O\}$, 
\item[{\rm (2)}] for an underlying immersion $k: \amalg S^1 \to \text{\bf R}^2 -\{O\}$ 
of $K$, the composition 
$\pi \circ k: \amalg S^1 \to S^1$ is an orientation preserving 
covering map of degree $m$ (where 
$\amalg S^1$ is the disjoint union of $\mu(K)$ circles), and 
\item[{\rm (3)}] $\pi|_{V(K)}: V(K) \to S^1$ is injective. 
\end{itemize}
A point $\theta$ of $S^1$ is called a {\it regular value\/} 
if $V(K) \cap \pi^{-1}(\theta) = \emptyset$.   By cutting 
$K$ along the half-line $\pi^{-1}(\theta)$ for a regular value $\theta$, 
we obtain a virtual braid diagram whose 
closure is $K$.  Such a virtual braid is unique up to conjugation.

{\bf Braiding Process {\rm (Proof of Proposition~3.1)}.} 
Let $K$ be a virtual link diagram and let $N_1, \dots, N_n$ be 
regular neighborhoods of the real crossings of $K$.  
By an isotopy of $\text{\bf R}^2$, we may assume that 
all $N_i$ ($i=1, \dots, n$) are in $\text{\bf R}^2 -O$, 
$\pi(N_i) \cap \pi(N_j) = \emptyset$ for $i \neq j$ and 
the restriction of $K$ to $N_i$ 
satisfies the condition of a braided virtual link diagram.  
Replace the remainder $K|_{W(K)}$ arbitrarily 
such that the result is a braided virtual link diagram. 
By Lemma~4.1,  $K$ is equivalent to this diagram.  
 
\section*{5. Proof of Theorem 3.2 and Proposition 3.4}

The terminologies ``braid moves'', ``right stabilizations'' and 
``right/left virtual exchange moves'' are also used 
for braided virtual link diagrams. 
(Conjugations are just braid moves.)  
These moves and their inverse moves are 
also called ${\rm VM0}$-, ${\rm VM2}$- and ${\rm VM3}$-moves.  
For example, the moves illustrated in Figure~9 are 
right stabilizations (${\rm VM2}$-moves) for braided virtual link  
diagrams.  
If two braided virtual link diagrams are related 
by a finite sequence of ${\rm VM0}$- and ${\rm VM2}$-moves, then 
we say that they are {\it virtually Markov equivalent in the strict sense\/}.  
If they are related 
by a finite sequence of ${\rm VM0}$-, ${\rm VM2}$- and ${\rm VM3}$-moves, then 
we say that they are {\it virtually Markov equivalent\/}.  

\fgVStabB

\noindent{\bf Lemma 5.1} {\it 
Let $K$ and $K'$ be braided virtual link diagrams 
(possibly of distinct degrees) 
such that $K'$ is obtained from $K$ by replacing $K|_{W(K)}$.  Then 
$K$ and $K'$ are virtually Markov equivalent in the strict sense.    
}

{\it Proof.} 
Let $N_1, \dots, N_n$ be regular neighborhoods of 
the real crossings of $K$ 
(and $K'$) and 
$W = W(K) = {\rm Cl}(\text{\bf R}^2 - \cup_{i=1}^n N_i)$.  
Take a common regular value $\theta_0 \in S^1$  for 
$K$ and $K'$ such that $\theta_0$ is not in $\pi(\cup_{i=1}^n N_i)$.    
Assume that there exists an arc/loop $a_i$ 
of $K|_W$ and the corresponding one $a'_i$ of $K'|_W$ 
such that $\sharp(a_i \cap \pi^{-1}(\theta_0)) \neq 
\sharp(a'_i \cap \pi^{-1}(\theta_0))$.    
Move a small part of $a_i$ or $a'_i$ toward the origin 
by a series of ${\rm VM0}$-moves corresponding 
to $\tau_i^2 =1$ and apply some ${\rm VM2}$-moves of virtual type 
so that $\sharp(a_i \cap \pi^{-1}(\theta_0)) = 
\sharp(a'_i \cap \pi^{-1}(\theta_0))$.
Thus we may assume that $\sharp(a_i \cap \pi^{-1}(\theta_0)) = 
\sharp(a'_i \cap \pi^{-1}(\theta_0))$ for every 
arc/loop $a_i$ 
of $K|_W$.
Let $k$ and $k'$ be underlying immersions 
$\amalg S^1 \to \text{\bf R}^2 -\{O\}$ of $K$ and $K'$ 
such that they are identical near the preimages of the real crossings.  
Let $I_1,\dots, I_s$ be intervals or circles in 
$\amalg S^1$  with $k(I_i)=a_i$ for $i=1,\dots, s$, and 
put $k_i=k|_{I_i}$.    Let $k'_1, \dots, k'_s$ be 
such immersions obtained from $K'$.  
Note that $\pi \circ k_i : I_i \to S^1$ and 
$\pi \circ k'_i : I_i \to S^1$ are orientation 
preserving immersions and $\pi \circ k_i|_{\partial I_i} = 
\pi \circ k'_i|_{\partial I_i}$.  
Since $a_i$ and $a'_i$ have the 
same degree with respect to $\theta_0$, we have a homotopy 
$\{k_i^s: I_i \to  \text{\bf R}^2 -\{O\} \}_{s\in [0,1]}$ 
between $k_i=k_i^0$ and $k'_i=k_i^1$ relative to the boundary 
$\partial I_i$ 
such that for each $s \in [0,1]$, 
$\pi \circ k_i^s : I_i \to S^1$ is an immersion.  
Taking such a homotopy generically with respect to the 
other arcs/loops of $K|_W$ (and $K'|_W$) and the 2-disks 
$N_1, \dots, N_n$, we have a finite sequence of 
${\rm VM0}$-moves transforming $a_i$ to $a'_i$ (recall the proof of Lemma~4.1).  
Applying this procedure inductively, we see that 
$K$ is transformed into $K'$ by ${\rm VM0}$-moves.  
\qed

Applying the above argument, we obtain Proposition~3.4.  

{\it Proof of Proposition~3.4.}  
In Figure~\ref{fig:VBraidD}, we show a process that 
$b \in VB_m$ is transformed into $\iota_1^0(b) \tau_1 \in VB_{m+1}$   
(the figure is for the case of $m=3$). 
The step (2) $\to$ (3) is a ${\rm VM2}$-move, up to 
${\rm VM1}$-moves.  
The other steps are 
${\rm VM0}$-moves and ${\rm VM1}$-moves. 
\qed

\fgVBraidD

\noindent{\bf Lemma 5.2} {\it 
Two braided virtual link diagrams with the 
same Gauss data are virtually Markov equivalent in the strict sense. 
}

{\it Proof.}  
Let $K$ and $K'$ be braided virtual link diagrams with the 
same Gauss data.  Let $N_1, \dots, N_n$ be regular neighborhoods 
(as in Figure~\ref{fig:Crossings}) of 
the real crossings $v_1, \dots, v_n$ of $K$, and 
$N'_1, \dots, N'_n$ be regular neighborhoods 
of the corresponding real crossings $v'_1, \dots, v'_n$ of $K'$.    

(Case 1)  Suppose that $\pi(N_1), \dots, \pi(N_n)$ and 
$\pi(N'_1), \dots, \pi(N'_n)$ appear in $S^1$ 
in the same (cyclic) order.   By an isotopy of $\text{\bf R}^2$,  
deform $K$ keeping the condition of a braided virtual link diagram 
such  that 
$N_1 = N'_1, \dots, N_n = N'_n$ and the restrictions of $K$ and $K'$ 
to these disks are identical.  
By Lemma~5.1, $K$ and $K'$ are 
virtually Markov equivalent in the strict sense.  

(Case 2)  Suppose that $\pi(N_1), \dots, \pi(N_n)$ and 
$\pi(N'_1), \dots, \pi(N'_n)$ do not appear in $S^1$ 
in the same (cyclic) order. 
It is sufficient to consider a special case that 
$\pi(N_1), \dots, \pi(N_n)$ and 
$\pi(N'_1), \dots, \pi(N'_n)$ appear in $S^1$ 
in the same order except a pair, say  
$\pi(N_1)$ and $\pi(N_2)$.  
Applying ${\rm VM0}$-moves, we may assume that $K$ is the closure of 
a virtual braid diagram which looks like 
the left one of Figure~\ref{fig:VBraidA}, where $b_1$ is a virtual 
braid diagram without real crossings and $b_2$ is a virtual 
braid diagram.  
The middle of 
the figure is obtained from the left by ${\rm VM0}$- and ${\rm VM2}$-moves.  
The right one is obtained from the middle by ${\rm VM0}$-moves.  By Case~1, 
the right one and $K'$ are virtually Markov equivalent 
in the strict sense.  Thus $K$ and $K'$ are 
virtually Markov equivalent in the strict sense.  
\qed
 
\fgVBraidA

Since the braiding process (given in \S~4) does not 
change the Gauss data of a 
virtual link diagram, we have the following.  

\noindent{\bf Corollary 5.3} {\it 
For a virtual link diagram $K$, a braided virtual link diagram 
obtained by the braiding process is unique up to 
virtual Markov equivalence in the strict sense. 
}

\fgVRMoveOriC

{\it Proof of Theorem~3.2.}  
The if part is obvious. We prove the only if part.  
Let $K$ and $K'$ be braided virtual link diagrams which are 
equivalent as virtual links.  
There is a finite sequence of virtual link diagrams from $K$ to $K'$ 
each step of which is one of the moves in Figure~\ref{fig:VRMoveOriC}.  
(For RI-moves and VI-moves, there are other cases of 
orientations of the arcs.  These cases are obtained from 
the moves in the figure by RII- and VII-moves. This is 
called {\it the Whitney trick\/}.  
For RIII-, VIII- and MI-moves, there are other cases of 
orientations on the arcs.  These cases are also obtained from 
the moves in the figure by RII- and VII-moves.)  
By use of VII-moves, an RIIc-move and an RIId-move 
are obtained from 
an Xa-move and an Xb-move in Figure~\ref{fig:VRMoveOriD}, respectively.  
Therefore,   
there is a finite sequence of virtual link diagrams 
$K= K_0, K_1, \dots, K_s =K'$ such that each $K_i$ 
is obtained from $K_{i-1}$ by an  
RIa-, RIb-, VI-, RIIa-, RIIb-, Xa-, Xb-, VIIa-, VIIb-, VIIc-, 
RIII-, VIII- or MI-move.  

\fgVRMoveOriD

Apply the braiding process to each $K_i$ and let $\widetilde{K_i}$ 
be a braided virtual link diagram with the same Gauss data 
as $K_i$.  Note that $\widetilde{K}_i$ is uniquely determined up to 
virtual Markov equivalence in the strict sense (Lemma~5.2).    
We assume that $\widetilde{K}_0 = K_0 =K$ and 
$\widetilde{K}_s = K_s =K'$.   
Then it is sufficient to prove that for each $i$ $(i=1, \dots, s)$, 
$\widetilde{K}_i$ and $\widetilde{K}_{i-1}$ are virtually 
Markov equivalent.  

If $K_i$ is obtained from $K_{i-1}$ by a VI-, VIIa-, VIIb-, VIIc-, 
 VIII- or MI-move, then $K_i$ and $K_{i-1}$ have the same Gauss data 
and so do $\widetilde{K}_i$ and $\widetilde{K}_{i-1}$.  
By Lemma~5.2,  $\widetilde{K}_i$ and $\widetilde{K}_{i-1}$ 
are virtually Markov equivalent.  

Suppose that $K_i$ is obtained from $K_{i-1}$ by an  
RIa-, RIb-, RIIa-, RIIb-, Xa-, Xb-, or  
RIII-move. 
Let $\Delta$ be a 2-disk in $\text{\bf R}^2$ where 
the move is applied, and let $\Delta^c$ be the complement 
of $\Delta$ in   $\text{\bf R}^2$ so that 
$K_i \cap \Delta^c = K_{i-1} \cap \Delta^c$.  

If the move is not an Xb-move, then we can 
deform $K_i$ and $K_{i-1}$ by an isotopy of $\text{\bf R}^2$ 
such that $K_i \cap \Delta$ and $K_{i-1} \cap \Delta$ 
satisfy the condition of a braided virtual link diagram.  
Apply the braiding process to the remainder 
$K_i \cap \Delta^c = K_{i-1} \cap \Delta^c$, and we have 
braided virtual link diagrams, 
say $\widetilde{K}'_i$ and $\widetilde{K}'_{i-1}$ 
such that 
$\widetilde{K}'_i \cap \Delta = K_i \cap \Delta$, 
$\widetilde{K}'_{i-1} \cap \Delta = K_{i-1} \cap \Delta$, and 
$\widetilde{K}'_i \cap \Delta^c = \widetilde{K}'_{i-1} \cap \Delta^c$.  
If the move is an RIa-, RIb-, or Xa-move, then $\Delta$ 
contains the origin $O$ of $\text{\bf R}^2$ and 
$\widetilde{K}'_i$ and $\widetilde{K}'_{i-1}$ are 
related by a right stabilization  
of positive/negative type or a right virtual exchange move.  
If the move is an RIIa-, RIIb-, or  
RIII-move, then $\Delta$ is disjoint from $O$ and  
$\widetilde{K}'_i$ and $\widetilde{K}'_{i-1}$ are 
related by a ${\rm VM0}$-move.  
Since $\widetilde{K}'_i$ has the same Gauss data as $K_i$, 
it is virtually Markov equivalent to $\widetilde{K}_i$ 
by Lemma~5.2.  
Similarly $\widetilde{K}'_{i-1}$ is 
virtually Markov equivalent to $\widetilde{K}_{i-1}$.  
Therefore  $\widetilde{K}_i$ and $\widetilde{K}_{i-1}$ 
are virtually Markov equivalent.  

If the move is an Xb-move, then deform 
$K_i$ and $K_{i-1}$ by an isotopy of $\text{\bf R}^2$ 
such that they are the closures of the (virtual) tangles 
depicted as (A1) and (B1) in Figure~\ref{fig:VBraidB}, 
say $K'_i$ and $K'_{i-1}$, where $b_1$ and $b_2$ are 
virtual braid diagrams.     
(First deform $K_i \cap \Delta$ and $K_{i-1} \cap \Delta$ 
such that they 
are as in the thick boxes of (A1) and (B1).  Then apply the 
braiding process to the remainder.)   
Let $\widetilde{K}'_i$ and $\widetilde{K}'_{i-1}$ be 
the closures of the virtual braid diagrams  
depicted as (A2) and (B2) in the figure.  
Note that $\widetilde{K}'_i$ has the same Gauss data 
as $K'_i$ and hence as $K_i$.   
Thus $\widetilde{K}'_i$ is virtually Markov equivalent to 
$\widetilde{K}_i$ (Lemma~5.2).  Similarly 
$\widetilde{K}'_{i-1}$ is virtually Markov equivalent to 
$\widetilde{K}_{i-1}$.  
On the other hand, 
$\widetilde{K}'_i$ and $\widetilde{K}'_{i-1}$ are related 
by a left virtual exchange move.  
Therefore  $\widetilde{K}_i$ and $\widetilde{K}_{i-1}$ 
are virtually Markov equivalent.   
\qed

\fgVBraidB

\section*{6. Proof of Propositions 3.3 and 3.5} 

We fix a positive integer $N$ and an integer $\alpha$.  
Let $K$ be a 
virtual link diagram and let 
$N_1, \dots, N_n$ be regular neighborhoods of the 
real crossings of $K$.  Let $a_1, \dots, a_s$ be 
the arcs/loops of $K|_{W(K)}$ as before.  
A {\it state\/} $S$ is assignment of elements 
of $\{ 1,2,\dots, N \}$ to the 
arcs/loops $a_1, \dots, a_s$.  
A state $S$ is {\it admissible\/} if 
at each crossing $v$ of $K$ the labels 
around $v$ are one of Figure~\ref{fig:Vweights}.  
Then we give $v$ an element 
of $\text{\bf Z}[q, q^{-1}]$ indicated in the figure, 
which is denoted by $g(K, S; v)$.  
For an admissible state $S$ of $K$, let 
$$G(K, S) = \prod_{v \in V(K)} g(K, S; v). $$

\fgVweights

\fgVsmooth

For an admissible state $S$ of $K$, 
let $K^S$ be the virtual link diagram obtained from 
$K$ by changing the crossing points of $K$ as in Figure~\ref{fig:Vsmooth}.  
Then each component, say $c$, of $K^S$ inherits a unique element 
of $\{1,2,\dots, N \}$ from the state $S$, 
which we denote by $S(c)$.  
Let 
$$H(K, S) = \prod_c q^{2S(c)-N-1}, $$
where $c$ runs over all components of $K^S$.  

We denote by $w(K)$ the number of positive crossings minus the 
number of negative crossings of $K$.  
For a virtual link diagram $K$, we define $Q_{N, \alpha}(K)$ by 
$$ Q_{N, \alpha}(K) = (-q^N)^{w(K)} 
\sum_S  G(K,S) H(K,S) \quad \in \text{\bf Z}[q, q^{-1}], $$
where $S$ runs over all admissible states of $K$.  For a virtual 
braid diagram $b$, we define $ Q_{N, \alpha}(b)$ by 
$ Q_{N, \alpha}({\rm closure~of~}b)$.  

\noindent{\bf Lemma 6.1} {\it 
If virtual braid diagrams $b$ and $b'$ are related by 
${\rm VM0}$-, ${\rm VM1}$- and ${\rm VM2}$-moves, then $Q_{N, \alpha}(b)=Q_{N, \alpha}(b')$.  
(If braided virtual link diagrams $K$ and $K'$ are 
virtually Markov equivalent in the strict sense, 
then 
$Q_{N, \alpha}(K)=Q_{N, \alpha}(K')$.)
}  

{\it Proof.}  
By a standard argument of state models (cf. 
\cite{\rJ, \rJKS, \rKauD, \rKauE, \rTu}), it is directly 
checked that $Q_{N, \alpha}$ does not change under each move. 
Details are left to the reader. 
\qed

{\it Remark.} 
The function $Q_{N, \alpha}$ is a modification of the state model 
of the braid invariant $T_S(b)$ given by V.~Turaev  \cite{\rTu}.   
$G(K,S)$ and $H(K,S)$ correspond to $\prod (f)$ and $\int_D f$ 
in \cite{\rTu}.  We changed $\int_D f$  into $H(K,S)$ so that 
the function does not change under a VI-move in 
Figure~\ref{fig:VRMoveOriC}.   
This yields loss of invariability under RIIc-, RIId-moves, that  
helps us to prove Propositions~3.3 and 3.5. 

{\it Proof of Proposition~3.3.}  
Let $b_1= \tau_1 \sigma_1^{-1} \tau_2 \tau_1 \sigma_1 \tau_2$ 
and $b_2= \tau_1 \sigma_1^{-1} \sigma_2^{-1} \tau_1 \sigma_1 \sigma_2 \in VB_3$.  
They are related by a right virtual exchange move.  
By a direct calculation,  we have 
$Q_{2,0}(b_1)=0$ and 
$Q_{2,0}(b_2)= q^{-3} -q^{-1} -q +q^3$.  
Therefore, a right virtual exchange move is not a consequence of 
${\rm VM0}$-, ${\rm VM1}$- and ${\rm VM2}$-moves.  
Let $b_3= \sigma_2^{-1} \tau_1 \sigma_2^{-1} \tau_1 \in VB_3$ 
and $b_4= \sigma_2^{-1} \sigma_1^{-1} \sigma_2^{-1} \sigma_1 \in VB_3$ .  
Then
$Q_{2,0}(b_3)= q^{-7} -q^{-5} -q^{-3} +2 q^{-1} +q$ and 
$Q_{2,0}(b_4)= q^{-1} +q$.  
Therefore, a left virtual exchange move is not a consequence of 
${\rm VM0}$-, ${\rm VM1}$- and ${\rm VM2}$-moves.    
\qed

{\it Proof of Proposition~3.5.}  
By a direct calculation, we have 
$Q_{2,0}(\tau_1 \sigma_1^{-1} \in VB_2)= 1 -q^{-2}$, 
$Q_{2,0}(\tau_2 \sigma_2^{-1} \sigma_1 \in VB_3)= -1 +q^2$, and 
$Q_{2,0}(\tau_2 \sigma_2^{-1} \sigma_1^{-1} \in VB_3)= 1 + q^{-6} -2 q^{-4}$. 
Therefore, we have the result.   
\qed
 
\section*{7. Welded Knots and Their Braid Presentation}

In this section a virtual link diagram is referred to as a {\it welded 
link diagram\/}.  We call the local move illustrated in the left hand side 
of Figure~\ref{fig:WRMoveA} a {\it ${\rm W}$-move}.  
Two welded link diagrams are {\it equivalent\/} 
{\it as welded links\/} 
if they are 
related by a finite sequence of virtual Reidemeister moves 
and W-moves.  
The equivalence class is called a {\it welded link} 
or a {\it welded link type}.  
It is easily verified that the oriented W-move 
illustrated in the right of Figure~\ref{fig:WRMoveA} 
is sufficient to 
realize all possible orientations for a W-move 
up to oriented moves in Figure~\ref{fig:VRMoveOriC}.  
    
\fgWRMoveA

We refer to a virtual braid diagram as a {\it welded 
braid diagram\/}. Recall that the welded braid group $WB_m$ is the quotient 
of $VB_m$ by adding the relations  
$\tau_i \sigma_{i+1} \sigma_i = \sigma_{i+1} \sigma_i \tau_{i+1}$ 
($i=1,\dots,m-2$) 
corresponding to W-moves.  

\noindent{\bf Proposition 7.1} {\it 
Every welded link type is represented by 
the closure of a welded braid diagram.  
}

{\it Proof.} This is a direct consequence of Proposition~3.1. \qed

\noindent{\bf Theorem 7.2} {\it 
Two welded braids (or welded braid diagrams)  
have equivalent closures as welded links 
if and only if they are related by a finite sequence of the 
following moves $({\rm WM1})$ -- $({\rm WM2})$ 
(or $({\rm WM0})$ -- $({\rm WM2})$): 
\begin{itemize}
\setlength{\itemsep}{-3pt} 
\item[{\rm (WM0)}] a welded braid move (which is a move   
corresponding to a defining relation of the welded braid group), 
\item[{\rm (WM1)}] a conjugation in the welded braid group, 
\item[{\rm (WM2)}] a right stabilization of positive, negative or virtual type, 
and its inverse operation. 
\end{itemize} 
}

\noindent{\bf Lemma 7.3} {\it 
A left stabilization of positive, negative or 
virtual type is a consequence of 
${\rm WM0}$-, ${\rm WM1}$- and ${\rm WM2}$-moves.  
}

{\it Proof.}  
If it is of virtual type, then it follows from 
Proposition~3.4.  
If it is of positive/negative type, then 
replace the virtual crossings of (2) in Figure~\ref{fig:VBraidD} 
with real crossings so that the step (4) $\to$ (5) is 
allowed in the welded braid group.  
\qed

\noindent{\bf Lemma 7.4} {\it 
A right/left virtual exchange move is a consequence of 
${\rm WM0}$-, ${\rm WM1}$- and ${\rm WM2}$-moves.  
}

{\it Proof.} A right virtual exchange move is realized by 
${\rm WM0}$-, ${\rm WM1}$- and ${\rm WM2}$-moves 
as follows: 
$$\begin{array} {ll}
b_1 \sigma_{m}^{-1} b_2 \sigma_m  
& = b_1 \sigma_{m}^{-1} \tau_m \tau_m  b_2 \sigma_m 
   \in WB_{m+1} \\
& \leftrightarrow b_1 \sigma_{m}^{-1} \tau_m \tau_{m+1} \tau_m  b_2 \sigma_m 
   \in WB_{m+2} \quad ({\rm WM1} + {\rm WM2})\\
& = b_1 \sigma_{m}^{-1} \tau_{m+1} \tau_m \tau_{m+1} b_2 \sigma_m 
   \in WB_{m+2} \\
& = b_1 \tau_{m+1} \tau_m \sigma_{m+1}^{-1} \tau_{m+1} b_2 \sigma_m 
   \in WB_{m+2} \\
& = \tau_{m+1} b_1  \tau_m  b_2 \sigma_{m+1}^{-1} \tau_{m+1} \sigma_m 
   \in WB_{m+2} \\
& \leftrightarrow b_1  \tau_m  b_2 \sigma_{m+1}^{-1} \tau_{m+1} \sigma_m \tau_{m+1} 
   \in WB_{m+2} \quad ({\rm WM1})  \\
& = b_1  \tau_m  b_2 \sigma_{m+1}^{-1} \tau_{m} \sigma_{m+1} \tau_{m} 
   \in WB_{m+2}  \\
& = b_1  \tau_m  b_2 \sigma_{m} \tau_{m+1} \sigma_{m}^{-1} \tau_{m} 
   \in WB_{m+2}  \\
& \leftrightarrow b_1  \tau_m  b_2 \sigma_{m} \sigma_{m}^{-1} \tau_{m} 
   \in WB_{m+1} \quad ({\rm WM1} + {\rm WM2}) \\
& = b_1  \tau_m  b_2 \tau_{m} \in WB_{m+1},  
\end{array}
$$  
where $b_1, b_2 \in WB_{m}$ (and we also denote by $b_i$ ($i=1,2$) the natural images 
$\iota_0^1(b_i) \in WB_{m+1}$ and $\iota_0^2(b_i) \in WB_{m+2}$). 
Similarly, a left virtual exchange move is realized by 
WM0-, WM1-moves and left stabilizations.  By Lemma~7.3, we have the result. \qed

We call a braided virtual link diagram a {\it braided welded link diagram}.  
Two  braided welded link diagrams are {\it welded Markov equivalent} if 
they are related by WM0- and WM2-moves.  
(WM1-moves are regarded as WM0-moves.) 
By Lemma~7.4, if two braided welded link diagrams are virtually Markov equivalent, 
then they are welded Markov equivalent.  

{\it Proof of Theorem~7.2.}  
The if part is obvious. We prove the only if part.  
Let $K$ and $K'$ be braided welded link diagrams which are 
equivalent as welded links.  
There is a finite sequence of welded link diagrams 
$K= K_0, K_1, \dots, K_s =K'$ such that each $K_i$ 
is obtained from $K_{i-1}$ by an  
RIa-, RIb-, VI-, RIIa-, RIIb-, Xa-, Xb-, VIIa-, VIIb-, VIIc-, 
RIII-, VIII-, MI- or W-move 
(in Figures~\ref{fig:VRMoveOriC}, \ref{fig:VRMoveOriD} 
and~\ref{fig:WRMoveA}).  
Apply the braiding process to each $K_i$ and let $\widetilde{K_i}$ 
be a braided welded link diagram with the same Gauss data 
as $K_i$.  By Lemma~5.2 (and Lemma~7.4),  
$\widetilde{K}_i$ is uniquely determined up to 
welded Markov equivalence.    
   We assume that $\widetilde{K}_0 = K_0 =K$ and 
$\widetilde{K}_s = K_s =K'$.   
It is sufficient to prove that for each $i$ $(i=1, \dots, s)$, 
$\widetilde{K}_i$ and $\widetilde{K}_{i-1}$ are welded   
Markov equivalent.  
In the proof of Theorem~3.2, we have already seen that 
$\widetilde{K}_i$ and $\widetilde{K}_{i-1}$ are welded   
Markov equivalent, except the case that  
$K_i$ is obtained from $K_{i-1}$ by a W-move.  
    Suppose that $K_i$ is obtained from $K_{i-1}$ by a W-move.  
Let $\Delta$ be a 2-disk in  $\text{\bf R}^2$ where the W-move 
is applied, and 
let $\Delta^c$ be the complement 
of $\Delta$ so that 
$K_i \cap \Delta^c = K_{i-1} \cap \Delta^c$.  
Deform $K_i$ and $K_{i-1}$ by an isotopy of $\text{\bf R}^2$ 
such that $K_i \cap \Delta$ and $K_{i-1} \cap \Delta$ 
satisfy the condition of a braided virtual (welded) link diagram.  
Apply the braiding process to the remainder 
$K_i \cap \Delta^c = K_{i-1} \cap \Delta^c$, and we have 
braided welded link diagrams, 
say $\widetilde{K}'_i$ and $\widetilde{K}'_{i-1}$ 
such that 
$\widetilde{K}'_i \cap \Delta = K_i \cap \Delta$, 
$\widetilde{K}'_{i-1} \cap \Delta = K_{i-1} \cap \Delta$, and 
$\widetilde{K}'_i \cap \Delta^c = \widetilde{K}'_{i-1} \cap \Delta^c$.  
$\widetilde{K}'_i$ and $\widetilde{K}'_{i-1}$ are 
related by a WM0-move 
corresponding to 
$\tau_k \sigma_{k+1} \sigma_k = \sigma_{k+1} \sigma_k \tau_{k+1}$.  
Since $\widetilde{K}'_i$ has the same Gauss data as $K_i$, 
it is welded Markov equivalent to $\widetilde{K}_i$.   
Similarly $\widetilde{K}'_{i-1}$ is 
welded Markov equivalent to $\widetilde{K}_{i-1}$.  
Therefore  $\widetilde{K}_i$ and $\widetilde{K}_{i-1}$ 
are welded Markov equivalent.  
\qed

{\it Remark.} 
(1) 
S.~Satoh \cite{\rSa} showed that welded links 
are related with ribbon surfaces in 4-space 
whose components are tori.  
From the point of view of \cite{\rSa}, welded braids 
are related with the motion group of a 
trivial link in 3-space 
(cf. \cite{\rGoldA, \rGoldB, \rMUY}).  

(2) 
When we use a move illustrated 
in Figure~\ref{fig:WRMoveB}, called a 
{\it ${\rm W}^*$-move\/}, instead of 
a W-move, we have another notion which is similar to 
a welded link.  
Define a group $WB^*_m$ by the quotient of 
$VB_m$ by the relations 
$\tau_i \sigma_{i+1}^{-1} \sigma_i^{-1} 
= \sigma_{i+1}^{-1} \sigma_i^{-1} \tau_{i+1}$ 
($i=1,\dots,m-2$),  
instead of 
$\tau_i \sigma_{i+1} \sigma_i 
= \sigma_{i+1} \sigma_i \tau_{i+1}$.  
Then we have results similar to those in this section. 
It should be noticed that we cannot use 
both of W-moves and ${\rm W}^*$-moves simultaneously.  
If we use both moves, every virtual (or welded) knot 
diagram changes into the unknot, 
\cite{\rGPV, \rKane, \rNelson}. 

\fgWRMoveB

\vskip0.3cm

\noindent {\it Address}: Department of Mathematics, Osaka City University, 
Sumiyoshi-ku, Osaka 558-8585, Japan  


\noindent {\it Current address} (until September 30, 2000): 
Department of Mathematics and Statistics, University of South Alabama, 
Mobile, AL 36688, USA  


\end{document}